\documentclass[journal]{IEEEtran}

\usepackage{booktabs} 
\usepackage{multirow}

\usepackage{tabularx}

\usepackage{makecell}

\usepackage{graphicx}

\usepackage{enumitem}

\usepackage{smartdiagram}
\usesmartdiagramlibrary{additions}

\graphicspath{{Figs/}}

\usepackage{pdfpages}
\usepackage{pdflscape}
\usepackage{amsmath}
\usepackage{amsthm}
\usepackage{amssymb}

\DeclareMathOperator*{\argmin}{arg\,min}

\usepackage{amsfonts}

\usepackage{color}

\usepackage{setspace}
\usepackage{enumitem}
\usepackage{xcolor}

\usepackage{bm}
\usepackage{bbm}

\usepackage{tikz}
\usepackage{pgf}
\usetikzlibrary{arrows,automata}
\usetikzlibrary{shapes}
\tikzstyle{data44}=[rectangle split,rectangle split parts=2,draw,text centered]

\usepackage{subcaption}

\DeclareMathAlphabet{\mathcal}{OMS}{cmsy}{m}{n}

\usetikzlibrary{matrix}

\tikzset{
  BarreStyle/.style = {opacity=.3,line width=14 mm,color=#1},
  node style ge/.style={},
  node style sp/.style={},
  yl/.style={},
  arrow style mul/.style={},
}

\newtheorem{remark}{Remark}

\usepackage{threeparttable,booktabs}

\usepackage[colorlinks,
            linkcolor=blue,
            anchorcolor=blue,
            citecolor=blue
            ]{hyperref}

\usepackage{cite}

\usepackage{mdwlist}

\usepackage{comment}

\usepackage{stmaryrd}

\usepackage[linesnumbered,ruled,vlined]{algorithm2e}

\SetKwInput{KwInput}{Input}                
\SetKwInput{KwOutput}{Output}

\let\oldnl\nl
\newcommand{\nonl}{\renewcommand{\nl}{\let\nl\oldnl}}

\AtBeginDocument{%
   \setlength\abovedisplayskip{4pt}
   \setlength\belowdisplayskip{4pt}}

\newcommand\smallmath[2]{#1{\raisebox{\dimexpr \fontdimen 22 \textfont 2
      - \fontdimen 22 \scriptscriptfont 2 \relax}{$\scriptscriptstyle #2$}}} 
\newcommand\smallotimes{\smallmath\mathbin\otimes}

\usepackage{stfloats}
\usepackage{tabularx}

\usepackage{textcomp}

\newcolumntype{P}[1]{>{\centering\arraybackslash}p{0.66cm}}
\usepackage{colortbl}

\usepackage{ragged2e}

\usepackage{changepage}

\newcommand{\hig}[1]{\textcolor{red}{#1}}

\usepackage{ragged2e}

\begin{document}

\title{Optimal Transmission Switching with Uncertainties from both Renewable Energy $\!\!\!$ and N-k Contingencies}



\author{Tong Han,~\IEEEmembership{}
        David J. Hill,~\IEEEmembership{Life Fellow, IEEE},
        and Yue Song,~\IEEEmembership{Member, IEEE}
\thanks{All authors are with the Department of Electrical and Electronic Engineering, University of Hong Kong, Hong Kong (e-mail: hantong@eee.hku.hk; dhill@eee.hku.hk; yuesong@eee.hku.hk).}
\thanks{D. J. Hill is also with the Department of Electrical and Computer Systems Engineering, Monash University, Melbourne, Australia (e-mail: david.hill3@monash.edu).}

}

\maketitle

\begin{abstract}
    This paper focuses on the $N-k$ security-constrained optimal transmission switching (OTS) problem for variable renewable energy (VRE) penetrated power grids. A new three-stage stochastic and distributionally robust OTS model is proposed. The first stage has the primary purpose to schedule the power generation and network topology based on the forecast of VRE. The second stage controls the power generation and voltage magnitudes of voltage-controlled buses in response to VRE uncertainty, and the third stage reacts to $N\!-\!k$ contingencies additionally by line switching and load shedding. 
    The VRE and $N\!-\!k$ contingencies, considering different availability of their probability distributions, are tackled by stochastic and distributionally robust optimization, respectively. By adopting stage-wise realization of uncertainties in VRE and contingencies, the associated corrective controls with different mechanisms can be handled separately and properly, which makes the proposed OTS model more realistic than existing two-stage ones. For solving the proposed OTS model, its tractable reformulation is derived, and a solution approach that combines the nested column-and-constraint generation algorithm and Dantzig–Wolfe procedure is developed. Finally, case studies include a simple IEEE network for illustrative purposes and then real system networks to demonstrate the efficacy of the proposed approach.
\end{abstract}


\IEEEpeerreviewmaketitle

\section*{Notation and Nomenclature}


$\!\!\!\!\!\!\!\!$
\textit{Notation}:$\!$ 
    \textit{(1)} For a vector $\bm{x}$, $\dot{\bm{x}}$, $\ddot{\bm{x}}$ and $\dddot{\bm{x}}$ are the vectors with the same meaning as vector $\bm{x}$ but for the first, second and third stage, respectively. Note that superscripts ``$\cdot$'', ``$\cdot\cdot$'' and ``$\cdot\!\!\cdot\!\!\cdot$'' are typical notations for time derivatives, but there can be no confusion here since the paper does not involve any derivatives. 
    \textit{(2)} For two vectors or matrices $\bm{x}$ and $\bm{y}$ with the same dimensions, $\bm{x}^2$ is the entry-wise square of $\bm{x}$, $\bm{x} \!\circ\! \bm{y}$ is the entry-wise product between $\bm{x}$ and $\bm{y}$.
    \textit{(3)} For a vector $\bm{x}$, $\Vert \bm{x} \Vert_{-2\infty}$ is the second smallest absolute value of entries of $\bm{x}$.

\subsection{Sets and Graphs}
\addcontentsline{toc}{section}{Nomenclature}
\begin{IEEEdescription}[\IEEEusemathlabelsep\IEEEsetlabelwidth{$V_1,V$}]
\item[$\!\!\!\!\!\!$$\mathcal{S}$] Set of sets of all identical lines between two buses.
\item[$\!\!\!\!\!\!$$\mathcal{U}_{\rm c}$, $\mathcal{U}_{\rm v}$] Set of conventional/VRE-based generators.
\item[$\!\!\!\!\!\!$$\mathcal{V}$, $\mathcal{E}$] Set of buses/branches. Multiple branches between the same pair of buses are handled separately. 
\item[$\!\!\!\!\!\!$$\mathcal{G}(\mathcal{V}, \mathcal{E}, r)$] The undirected multigraph representing the transmission network topologically, and $r(e)$ with $e\!\in\! \mathcal{E}$ determines the two buses linked by branch $e$. 
\end{IEEEdescription}

\subsection{Scalars}
\addcontentsline{toc}{section}{Nomenclature}
\begin{IEEEdescription}[\IEEEusemathlabelsep\IEEEsetlabelwidth{$V_1,V$}]
\item[$\!\!\!\!\!\!$$k_{\rm max}$] The maximal number of fault components.
\item[$\!\!\!\!\!\!$$n_{\rm f}^i$] Number of sample points for the generation cost function of generator $i$.
\item[$\!\!\!\!\!\!$$n_{\rm p}$] Number of segments of the cosine approximation.
\item[$\!\!\!\!\!\!$$n_{\rm s}$] Number of segments of the approximation to output power constraints of VRE-based generators.
\item[$\!\!\!\!\!\!$$\Delta t$] Duration of corrective control.
\item[$\!\!\!\!\!\!$$M, N$]  Big-M constant, total number of components.
\end{IEEEdescription}

\subsection{Vectors}
\addcontentsline{toc}{section}{Nomenclature}
\begin{IEEEdescription}[\IEEEusemathlabelsep\IEEEsetlabelwidth{$V_1,V$}]
\item[$\!\!\!\!\!\!$$\bm{0}_{n}$, $\bm{1}_{n}$] The $n$-dimensional vector of zeros/ones, and ``$n$" is omitted if it is determinable by involved operations.

\item[$\!\!\!\!\!\!$$\bm{o} \!\in\! \mathbb{B}^{N} $]  Parameterization of $N-k$ contingencies. Entry values of 1/0 indicate normal/failure states of associated components.
\item[$\!\!\!\!\!\!$$\bm{o}_{\rm c}$, $\!\!\bm{o}_{\rm b}$, $\!\!\bm{o}_{\rm v}$]  Sub-vectors of $\bm{o}$ associated with $\mathcal{U}_{\rm c}$, $\mathcal{E}$ and $\mathcal{U}_{\rm v}$.

\item[$\!\!\!\!\!\!$$\bm{p}_{\rm c} \!\!+\!\! j \bm{q}_{\rm c}$] Complex power outputs of conventional generators.
\item[$\!\!\!\!\!\!$$\bm{p}_{\rm c_+}$, $\bm{p}_{\rm v_+}$] Upward regulations of active power outputs of contentional/VRE-based generators.
\item[$\!\!\!\!\!\!$$\bm{p}_{\rm c_-}$, $\bm{p}_{\rm v_-}$] Downward regulations of active power outputs of contentional/VRE-based generators.  
\item[$\!\!\!\!\!\!$$\tilde{\bm{p}}_{\rm c}^i$] Sample points within the range of active power outputs of concerns for generator $i$.
\item[$\!\!\!\!\!\!$$\bm{p}_{\rm c}^{\rm max}\!$, $\!\bm{p}_{\rm c}^{\rm min}$] Max/Min active power outputs of generators in $\mathcal{U}_{\rm c}$.  
\item[$\!\!\!\!\!\!$$\bm{p}_{\rm d} \!\!+\!\! j \bm{q}_{\rm d}$] Complex load powers.
\item[$\!\!\!\!\!\!$$\bm{p}_{\rm d_{\Delta}}$, $\bm{p}_{\rm d_{\Delta}}^{\rm max}$] Amounts of load shedding and their upper bounds.
\item[$\!\!\!\!\!\!$$\bm{p}_{\rm fb} \!\!+\!\! j \bm{q}_{\rm fb}$] Complex powers at the starting buses of branches.
\item[$\!\!\!\!\!\!$$\bm{p}_{\rm tb} \!\!+\!\! j \bm{q}_{\rm tb}$] Complex powers at the end buses of branches.
\item[$\!\!\!\!\!\!$$\bm{p}_{\rm v} \!\!+\!\! j \bm{q}_{\rm v}$] Complex power outputs of VRE-based generators. 
\item[$\!\!\!\!\!\!$$\bm{p}_{\rm v}^{\rm max}\!$, $\!\!\bm{s}_{\rm v}^{\rm max}$] Maximal active power outputs and MVA capacities of VRE-based generators, respectively.

\item[$\!\!\!\!\!\!$$\bm{q}_{\rm c}^{\rm max}\!$, $\!\bm{q}_{\rm c}^{\rm min}$] Max/Min reactive power outputs of $\mathcal{U}_{\rm c}$. 
\item[$\!\!\!\!\!\!$${\scriptstyle \Delta} \bm{q}_{\rm fb}$, $\!{\scriptstyle\Delta} \bm{q}_{\rm tb}$] Components of $\bm{q}_{\rm fb}$ and $\bm{q}_{\rm tb}$ which are coupled with changes in voltage magnitudes, respectively.

\item[$\!\!\!\!\!\!$${\bm{r}}_{\rm c_+}$, ${\bm{r}}_{\rm c_-}$] Upward/downward ramp rate limits of $\mathcal{U}_{\rm c}$.
\item[$\!\!\!\!\!\!$${\bm{r}}_{\rm v_+}$, ${\bm{r}}_{\rm v_-}$] Upward/downward ramp rate limits of $\mathcal{U}_{\rm v}$.

\item[$\!\!\!\!\!\!$$\bm{s}_{\rm b}^{\rm max}$] Power capacities of branches. 

\item[$\!\!\!\!\!\!$$\bm{v} \angle \bm{\theta}$] Complex voltages of buses. 
\item[$\!\!\!\!\!\!$$\bm{v}_{\rm c}$] Voltage magnitudes of voltage-controlled buses. 
\item[$\!\!\!\!\!\!$${\bm{v}}_{\rm c_{\Delta}}$, ${\bm{q}}_{\rm v_{\Delta}}$] Regulations of $\bm{v}_{\rm c}$ and $\bm{q}_{\rm v}$, respectively. 
\item[$\!\!\!\!\!\!$$\bm{v}_{\rm max}$, $\!\!\bm{v}_{\rm min}$] Max/Min voltage magnitudes of buses.

\item[$\!\!\!\!\!\!$$\bm{w}_{\rm c_+}$, $\bm{w}_{\rm c_-}$] Upward/downward regulation costs of $\mathcal{U}_{\rm c}$.
\item[$\!\!\!\!\!\!$$\bm{w}_{\rm d}$] Load shedding costs. 
\item[$\!\!\!\!\!\!$$\bm{w}_{\rm v_+}$, $\bm{w}_{\rm v_-}$] Upward/downward regulation costs of $\mathcal{U}_{\rm v}$.

\item[$\!\!\!\!\!\!$$\bm{z} \in \mathbb{B}^{|\mathcal{E}|}$] Statuses of branches. With $\bm{z} \!=\! [\bm{z}_k]_{k \in \mathcal{E}} $, $\bm{z}_{k} \!=\! 1$ indicates line $k$ is switched on and $\bm{z}_k\!=\! 0$ otherwise.
\item[$\!\!\!\!\!\!$$\bm{z}_{+}$, $\bm{z}_{-}$]  Action signs of switching on/off lines. A value of 1/0 means a/no switching action performed.

\item[$\!\!\!\!\!\!$$\bm{\varepsilon}$] Forecast errors of $\bm{p}_{\rm v}^{\rm max}$.
\item[$\!\!\!\!\!\!$$\bm{\alpha}$, $\bm{\vartheta}$, $\bm{\rho}$]  Electrical flow injections at nodes, vertex potentials, and electrical flows of edges all of the electrical flow network for network connectedness.    
\item[$\!\!\!\!\!\!$$\tilde{\bm{\eta}}^i$] Generation costs at the sample points of generator $i$. 
\item[$\!\!\!\!\!\!$$\bm{\theta}_{\rm max}$] Maximum phase angle differences of branches.
\item[$\!\!\!\!\!\!$$\bm{\phi}_{\rm v}^{\rm min}$] Minimum power factors of VRE-based generators. 
\item[$\!\!\!\!\!\!$$\bm{\varphi}_+$, $\bm{\varphi}_{-}$] Cosine approximation variables for branches, i.e., approximations to $\cos(\bm{E}^T \bm{\theta})$ and $\cos(- \bm{E}^T \bm{\theta})$. 
\end{IEEEdescription}

\subsection{Matrices}
\addcontentsline{toc}{section}{Nomenclature}
\begin{IEEEdescription}[\IEEEusemathlabelsep\IEEEsetlabelwidth{$V_1,V$}]
\item[$\!\!\!\!\!\!$$\bm{E}$, $\tilde{\bm{E}}$] Oriented incidence matrix of graph $\mathcal{G}$ with each branch assigned an arbitrary and fixed orientation, and its entry-wise absolute value.
\item[$\!\!\!\!\!\!$$\bm{E}_+$, $\!-\!\bm{E}_-$] Formed by replacing all $-1$/$1$ entries of $\bm{E}$ by $0$.  
\item[$\!\!\!\!\!\!$$\bm{E}_{\rm c}$] Incidence matrix between buses and conventional generators.
\item[$\!\!\!\!\!\!$$\bm{E}_{\rm d}$] Incidence matrix between buses and loads. 
\item[$\!\!\!\!\!\!$$\bm{E}_{\rm v}$] Incidence matrix between buses and VRE-based generators.  
\item[$\!\!\!\!\!\!$$\bm{G} \!+\!\! j \bm{B}$] Diagonal matrix formed by admittance of each branch.
\item[$\!\!\!\!\!\!$$\underline{\bm{G}} \!+\! j\underline{\bm{B}}$] Diagonal matrix formed by components of ground admittance of each bus contributed by bus shunts and transformers.  
\item[$\!\!\!\!\!\!$$\underline{\bm{G}}_{\rm b} \!+\! j \underline{\bm{B}}_{\rm b} $] Diagonal matrix formed by half ground admittance of each branch contributed by line charges.
\item[$\!\!\!\!\!\!$$\tilde{\bm{J}}$, $\bm{J}_{n(\!\times\! m)}$] Equal to matrix $\bm{J}$ in \cite{4-995-ea}, $n \times\! n$ ($\!\times m$) all-ones matrix. 
\item[$\!\!\!\!\!\!$$\bm{K}_{n \times m}$] The $n \times m$ matrix where all main diagonal entries are 1, all adjacent entries to the right of the main diagonal are $-1$, and the others are 0.
\end{IEEEdescription}

\section{Introduction}

\IEEEPARstart{M}{any} power grids are achieving increasingly high renewable penetration to mitigate global climate change. Variable renewable energy (VRE) is degrading the effectiveness of purely generation-side regulations due to the decreased generation dispatchablility. 
In this context, grid-side flexibility is expected to be leveraged more widely \cite{4-970}. 
Optimal transmission switching (OTS) is the primary means of leveraging the flexibility in network topology to improve economic efficiency of transmission networks. It can be generally seen as the augmentation of economic dispatch (ED) or unit commitment (UC) with optimization of network topology \cite{4-62, 4-670}. However, high-level VRE penetration brings new challenges for OTS.

Generally, it is indispensable to consider the $N\!-\!1$ or $N\!-\!k$ security criterion in OTS problems\cite{4-63, 4-670, 4-1140}. Meanwhile, for VRE penetrated power grids, the power outputs of VRE-based generators feature high intermittence and volatility and are hard to forecast accurately. Deterministic OTS models probably fail to ensure cost efficiency and operational security since they optimize with only one scenario. Thus for OTS with high VRE penetration, the awareness of VRE uncertainty is crucial to ensuring operational security in all possible scenarios of VRE outputs while providing statistically optimal decisions on generation dispatch and network topology. 
However, this simultaneous consideration of $N\!-\!k$ security and VRE uncertainty makes OTS problems much harder to tackle. For OTS formulation, we face at least the challenge of simultaneously and properly modeling the corrective control reacting to contingencies and VRE uncertainty. For computational tractability, $N\!-\!k$ security with numerous components of practical systems and VRE uncertainty cause double computation burdens.

Accordingly, in this paper, we devote to the OTS problem of VRE penetrated power grids, pursuing simultaneous, proper and efficient treatment of VRE uncertainty and $N\!-\!k$ security.

\subsection{Literature Review}

Most existing OTS models tackle VRE uncertainty under the framework of stochastic optimization (SO) or robust optimization (RO). In the stochastic OTS models, VRE outputs are modeled as uncertain parameters with known probability distributions (PDs). Ref. \cite{4-1137} developes a chance-constrained OTS model where the chance constraints are approximated by the sample average approach. The stochastic OTS models with VRE uncertainty handled by the point estimation method are proposed in \cite{4-79, 4-1223, 4-1158}. Two-stage stochastic OTS formulations are developed in \cite{4-1139} and \cite{4-1142} which, differing from the single-stage version \cite{4-1137, 4-79, 4-1223, 4-1158}, capture the corrective control reacting to VRE uncertainty. Unlike deterministic OTS \cite{4-62, 4-670}, the stochastic OTS can achieve optimal average performance. However, the computational challenge of SO models becomes significant when considering the $N-k$ security or for high-level VER penetrated power grids with numerous VRE-based generators. In the robust OTS models, VRE outputs are treated as uncertain parameters in an uncertainty set that ignores all distributional knowledge except for the support. Ref. \cite{4-1226} proposes an interval-based two-stage robust OTS model which maximizes the span of the possible wind variation interval. A two-stage robust OTS model for wind integrated hybrid AC/DC power grids is developed in \cite{4-1127}, where both uncertainties in wind generation and generator failure are considered. In contrast to stochastic OTS, the robust OTS models are computationally cheaper and the support is easily accessible. However, they often raise concerns of over-conservatism, i.e., pursuing better performance in the rare worst-case scenario while sacrificing the average performance. 
 
A common approach to consider $N\!-\!1$ security is to impose the operational feasibility on a pre-defined contingency set \cite{4-63, 4-670, 4-79, 4-1142, 4-1248}, which however, becomes computationally intensive for uncertainty-aware OTS of large-scale power grids. More importantly, this approach can be prohibitive for $N\!\!-\!\!k$ security due to massive possible contingencies. Another approach is to treat contingencies as uncertainties and formulate OTS as an RO model \cite{4-1127, 4-1140}. This approach is computationally efficient but can result in conservative solutions as the robust model to tackle VRE uncertainty. Recently, in UC and distribution network reconfiguration problems, $N\!-\!k$ security has been tackled under the framework of distributionally robust optimization (DRO) \cite{4-1195, 4-1214}, which prevents over-conservatism since the solution is induced by the worst-case distribution instead of the worst-case scenario in RO. However, there remains a lack of embedding $N\!-\!k$ security into OTS models under the framework of DRO to prevent conservatism while retaining computational tractability.

How to address the uncertainties in contingencies and VRE simultaneously also deserves a further study. All existing OTS models considering both these two kinds of uncertainties feature a two-stage stochastic/robust model \cite{4-1142, 4-1139, 4-1127}, with the implicit assumption that uncertainty realization of VRE and contingencies occurs simultaneously and their respective corrective control actions can be merged. In reality, however, these two uncertainties generally happen separately, and their corrective controls also have different mechanisms and thus being uncombinable. For example, corrective line switching control is adequate for reacting to contingencies but generally not for VRE uncertainty \cite{ 4-538, 4-779, 4-1248}. Generators engaging in different corrective controls may also vary. Moreover, under the framework of DRO, related works for OTS and other operational optimization problems tackle only one of uncertainties of VRE and contingencies \cite{4-1217, 4-1197, 4-1228, 4-1157, 4-1195, 4-1214}. 

The power flow equations employed by OTS models impact the practicality of OTS. The AC power flow model or its convex relaxations generally ensure AC feasibility and satisfactory solution optimality \cite{4-61, 4-1516, 4-1517}. 
These models, however, cause computational intractability for real-scale systems when both uncertainties in contingencies and VRE or even either is considered in OTS. 
Moreover, OTS with contingency uncertainty is basically some form of security-constrained optimal power flow (SCOPF) \cite{4-1518}. Some advanced and accurate power flow approximations, such as the sequential linearization \cite{4-1519} and successive linear programming approximation \cite{4-1520}, have been developed for common SCOPF problems. However, the usefulness of these approximations for OTS problems is highly impacted by the treatment of uncertainties. 
Therefore, most uncertainty-aware OTS models adopt less accurate yet more computationally tractable power flow approximations, including the DC power flow model \cite{4-79, 4-1248}, linearized power flow model \cite{4-1139}, and piecewise linear power flow model \cite{4-1142}. 
Two limitations of using these rough approximations are that the obtained solutions can be AC-infeasible, 
and the line switching actions can be far from the optimum in terms of the real AC model. 
For the first one, AC feasibility can be restored by common techniques before applying solutions to the real AC system \cite[Ch. 6]{4-490}. 
Due to the second limitation, the benefit from line switching evaluated based on the approximated model can be impacted or even become detriment for the real AC system. 
This impact of using the rough approximations, though neglected by most existing works, should be evaluated.

\subsection{Contributions}

In light of the challenges facing OTS and incompleteness of existing works, we proposes a new OTS model to simultaneously and properly treat VRE uncertainty and $N-k$ security, and a solution approach that enables utilization of parallel computing to improve computational efficiency. 
Specifically, the main contributions are summarized as follows: 

${\bullet}$ \textit{A three-stage stochastic and distributionally robust OTS (TSDR-OTS) model is developed.} The first stage has the primary purpose to schedule the power generation and network topology based on the forecast of VRE. 
The uncertainties of VRE and contingencies are considered at the second and third stage, respectively, such that their corrective controls can be modeled separately and properly. Corrective actions reacting to VRE uncertainty include generation and voltage regulation, and the third stage also leverages load shedding and corrective line switching reacting to $N\!\!-\!k$ contingencies. In the second stage, VRE uncertainty is addressed based on SO given that the PD of VRE forecast errors can be estimated relatively accurately based on available historical data. In contrast, an accurate PD of component failures is generally unattainable since the failure probabilities are close to zero and most contingencies have not actually happened in historical data. Therefore, at the third stage, contingency uncertainty is handled based on DRO. Consequently, the proposed TSDR-OTS model can achieve an optimal average performance regarding VRE uncertainty and optimal worst-case average performance regarding contingency uncertainty.

${\bullet}$ \textit{Tractable formulation of the TSDR-OTS model is derived for which an efficient solution approach, \hig{adapted from the enhanced column-and-constraint generation (CCG) algorithm in \cite{4-1213}}, is developed.} Since the TSDR-OTS model, with a complicated structure, cannot be solved directly, the Karush–Kuhn–Tucker (KKT) condition, strong duality and scenario-based approximation are employed to derive a traceable reformulation. The solution approach, combines the nested CCG algorithm and the Dantzig–Wolfe (DW) procedure and is more computationally efficient than alternatives. Particularly, compared with the CCG algorithm with a DW procedure developed in \cite{4-1213}, our approach \hig{additionally} contains an inner CCG loop, which utilizes strong duality to improve computational efficiency, to handle the more tricky subproblems in the DW procedure. \hig{These subproblems, due to the existence of binary decision variables, cannot be transformed into a tractable form handled by off-the-shelf solvers directly as in \cite{4-1213}, thus hindering adoption of the algorithm in \cite{4-1213}.}


\section{Mathematical Formulation}

\begin{figure}[h]
	\centering
 	\includegraphics[scale=0.8]{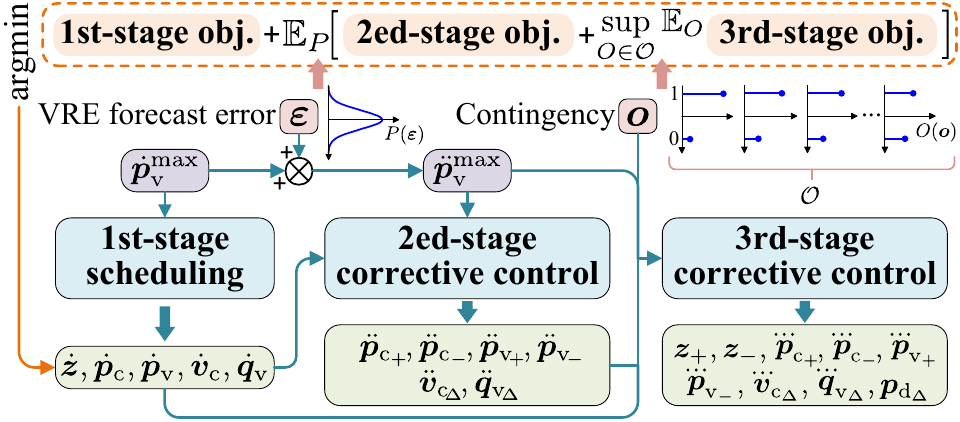}  
	\caption{Overall framework of the TSDR-OTS model.}
	\label{fig-3-6-3}
\end{figure}

This section formulates the TSDR-OTS model. Fig. \ref{fig-3-6-3} shows the overall framework of the TSDR-OTS model. 
It is noted that this work only focuses on the simpler case of OTS where UC variables are not considered, with the assumption that they are already determined. 
The first stage has the primary purpose to schedule the power generation and network topology based on the forecast of VRE $\dot{\bm{p}}_{\rm v}^{\max}$. Voltage magnitudes of voltage-controlled buses and reactive power outputs of VRE-based generators are optimized simultaneously. Then, the second stage finds the corrective controls of the power generation and voltage magnitudes of voltage-controlled buses, reacting to the VRE forecast errors $\bm{\varepsilon}$ due to VRE uncertainty. Corrective line switching is not considered since the main target of this stage is to correct the mild power imbalance ensuing from the VRE forecast errors. When no contingency occurs, i.e., $\bm{o} \!=\! \bm{1}$, the system operates with the topology scheduled at the first stage, and the power generation and voltage magnitudes corrected by the second stage. The third-stage corrective control is inactive in this case. When an $N\!-\!k$ contingency that can be simultaneous branch and/or generator failures, occurs, i.e, $\bm{o} \!\neq\! \bm{1}$, the third stage responds to it by corrective controls including line switching, generation redispatch, voltage regulation, and load shedding.
 
In the TSDR-OTS model, the VRE forecast errors $\bm{\varepsilon}$ are assumed to follow a PD denoted as $P$, which can be obtained from historical data; and contingency uncertainty $\bm{o}$ follows a PD denoted as $O$, which belongs to a certain ambiguity set $\mathcal{O}$. With these uncertainties and stage-specific objectives defined later, the first-stage scheduling aims to optimize the overall performance integrating the objectives of all three stages, with $\bm{\varepsilon}$ and $\bm{o}$ handled by SO and DRO respectively. 
Accordingly, the TSDR-OTS model yields the optimal operational solution which achieves overall optimality regarding normal operational performance with forecast VRE, average performance regarding VRE uncertainty, and worst-case average performance regarding contingency uncertainty.


\subsection{First-Stage Formulation}
 
\subsubsection{Objective} The first-stage objective is to minimize the sum of the individual costs of all conventional generators, i.e., 
\begin{equation}\label{eq-3-6-first-stage-obj-linear}
      \min~ \sum\nolimits_{ i \in \mathcal{U}_{\rm{c}} } \eta_{i} ( \dot{\bm{p}}_{{\rm{c}}, i} )  =  \sum\nolimits_{ i \in \mathcal{U}_{\rm{c}} } (\tilde{\bm{\eta}}^{i})^T \tilde{\bm{\lambda}}^{i} 
\end{equation} 
with the interpolation variables $\tilde{\bm{\lambda}}^i \!\in\! \mathbb{R}^{n_{\rm f}^i}$ satisfying 
\begin{equation}\label{eq-3-6-first-stage-obj-liearcons} 
    \dot{\bm{p}}_{{\rm c}, i} = (\tilde{\bm{p}}^i_{\rm c})^T \tilde{\bm{\lambda}}^i, \bm{0} \leq \tilde{\bm{\lambda}}^i \leq \bm{1}, \bm{1}^T \tilde{\bm{\lambda}}^i = 1 ~~ \forall i \in \mathcal{U}_{\rm c}
\end{equation}
where $\eta_i(\!\cdot)\!$ is the generation cost function of generator $i$, which is generally nonlinear and we assume that it can be approximated by a convex piecewise linear function; $\tilde{\bm{p}}_{\rm c}^i \in \mathbb{R}^{n_{\rm f}^i} $ and $\tilde{\bm{\eta}}^i \in \mathbb{R}^{n_{\rm f}^i}$ are known parameter vectors. Here the ``$\lambda$" approximation, one of the best-performing ones according to \cite{4-1240}, is adopted. It encodes a convex hull from the collection of sample points by introducing an interpolation variable in the range of $[0,1]$ for each sample point in the piecewise linear function and linking all of these interpolation variables together with the last constraint in (\ref{eq-3-6-first-stage-obj-liearcons}). The active power and generation cost of the interpolated point are recovered with the first constraint in (\ref{eq-3-6-first-stage-obj-liearcons}) and $(\tilde{\bm{\eta}}^{i})^T \tilde{\bm{\lambda}}^{i}$ in (\ref{eq-3-6-first-stage-obj-linear}), respectively.

\subsubsection{Power Flow Constraints} The AC power flow constraints are formulated as 
\begin{subequations}\label{eq-3-6-ac-power-1}
        \begin{align}
            &\!\!\!\!\!\!\!
            \begin{bmatrix}
                \! \dot{\bm{p}}_{\rm fb} \!\!\!  \\
                \! \dot{\bm{p}}_{\rm tb} \!\!\!  \\
                \! \dot{\bm{q}}_{\rm fb} \!\!\!  \\
                \! \dot{\bm{q}}_{\rm tb} \!\!\!
            \end{bmatrix}
            \!\!\!=\!\!\!\!
            \left(\!\!
            \begin{bmatrix}
                \!\bm{G}      \!\bm{E}_+^{T} \\
                \!\bm{G}      \!\bm{E}_-^{T}\!\! \\
                \!\!-\!\bm{B} \!\bm{E}_+^{T}\!\! \\
                \!\!-\!\bm{B} \!\bm{E}_-^{T}\!\!
            \end{bmatrix}
            \!\!\! \dot{\bm{v}}^{2}
            \!\!\!-\!\!\!
            \begin{bmatrix}
                \! \bm{G}      \!\!&\!\! \bm{B}     \!\!\!\\
                \! \bm{G}      \!\!&\!\!\!\!- \bm{B} \!\! \\
                \!\! -\bm{B} \!\!&\!\! \bm{G}     \!\!\!\\
                \!\! -\bm{B} \!\!&\!\!\!-\! \bm{G} \!\!
            \end{bmatrix}\!\!\!\!
            \begin{bmatrix}
                \!\cos(\! \bm{E}^{T} \!\dot{\bm{\theta}} \!) \!\! \\
                \!\sin(\! \bm{E}^{T} \!\dot{\bm{\theta}} \!) \!\!
            \end{bmatrix}
            \!\!\!\circ\!\!
            \bm{(}\!\bm{1}_4 \!\!\smallotimes\!\! 
            \psi\!(\! \dot{\bm{v}} \!) \!\bm{)}
            \!\!\! \right) 
            \!\!\!\circ\!\! (\! \bm{1}_{\!4} \!\!\smallotimes\!\! \dot{\bm{z}} )     \label{eq-3-6-ac-power-1:1} 
            \\
            &\!\!\!\!\!\!\! \begin{bmatrix}
                \!\! \bm{E}_{\rm{c}} \dot{\bm{p}}_{\rm{c}} \!\!+\!\! \bm{E}_{\!\rm{v}} \dot{\bm{p}}_{\rm{v}} \!\!-\!\! \bm{E}_{\rm d} \dot{\bm{p}}_{\rm{d}} \!\! \\
                \!\! \bm{E}_{\rm{c}} \dot{\bm{q}}_{\rm{c}} \!\!+\!\! \bm{E}_{\!\rm{v}} \dot{\bm{q}}_{\rm{v}} \!\!-\!\! \bm{E}_{\rm d} \dot{\bm{q}}_{\rm{d}} \!\!
            \end{bmatrix}
            \!\!\!-\!\!\!
            \begin{bmatrix}
                \!\! \bm{E}_{\!+} \dot{\bm{p}}_{\rm{fb}}  \!\!\!+\!\! \bm{E}_{\!-} \dot{\bm{p}}_{\rm{tb}} \!\!\! \\
                \!\! \bm{E}_{\!+}  \dot{\bm{q}}_{\rm{fb}} \!\!\!+\!\! \bm{E}_{\!-} \dot{\bm{q}}_{\rm{tb}} \!\!\!
            \end{bmatrix}
            \!\!\!=\!\!\! 
            \begin{bmatrix}
                \tilde{\bm{E}} \underline{\bm{G}}_{\rm b} \dot{\bm{z}} \!+\! \underline{\bm{G}} \bm{1}   \!\!\\
                \!\!-\! \tilde{\bm{E}} \underline{\bm{B}}_{\rm b} \dot{\bm{z}} \!\!-\!\! \underline{\bm{B}} \bm{1} \!\!
            \end{bmatrix}
            \!\!\!\circ\!\!\!
            \begin{bmatrix}
                \!\dot{\bm{v}}^{2}\!\! \\
                \!\dot{\bm{v}}^{2}\!\! 
            \end{bmatrix}\!\!\!\!
             \label{eq-3-6-ac-power-1:2}    
        \end{align}  
\end{subequations}
which is obtained by reformulating the AC power flow equations with voltages in polar coordinates to separate the power flow terms related to different admittance components, and considering the branch status parameterized by $\dot{\bm{z}}$. Here $\psi( \dot{\bm{v}} ) \!=\! (\bm{E}_+^T \dot{\bm{v}} ) \!\circ\! (\bm{E}_-^T  \dot{\bm{v}} )$; (\ref{eq-3-6-ac-power-1:1}) are branch power flow equations; for lines being switched on, i.e., with the corresponding entries in $\dot{\bm{z}}$ equal to 1, their branch power flow equations are valid, or otherwise, the branch and charging power are both forced to $0$; and (\ref{eq-3-6-ac-power-1:2}) represents power balance of each bus. 

For computational tractability, we linearize (\ref{eq-3-6-ac-power-1:1}) using the cold-start linear-programming approximation (LPA) \cite{4-1163} (with slight modification). Then, by further linearizing all bilinear terms w.r.t. $\bm{z}$ by the Big-M method, (\ref{eq-3-6-ac-power-1:1}) is approximated as
\begin{subequations}\label{eq-3-6-lpac}
        \begin{align}
            &  ( \dot{\bm{z}} \!-\! \bm{1} ) M  \! \leq \! \bm{G} \bm{1} \!\!-\!\! ( \bm{G} \dot{\bm{\varphi}}_{+} \!+\! \bm{B} \bm{E}^T\! \dot{\bm{\theta}} ) \!-\! \dot{\bm{p}}_{\rm fb} \!\leq\!  ( \bm{1} \!-\! \dot{\bm{z}} ) M   \label{eq-3-6-lpac:1}\\
            &  ( \dot{\bm{z}} \!-\! \bm{1} ) M  \! \leq \! \bm{G} \bm{1}  \!\!-\!\! ( \bm{G} \dot{\bm{\varphi}}_{+} \!-\! \bm{B} \bm{E}^T\! \dot{\bm{\theta}} ) \!-\! \dot{\bm{p}}_{\rm tb} \!\leq\!  ( \bm{1} \!-\! \dot{\bm{z}} ) M   \label{eq-3-6-lpac:2}\\
            &  (\! \dot{\bm{z}} \!\!-\!\! \bm{1} \!) M  \!\! \leq \!\! - \bm{B} \bm{1}  \!\!-\!\! (\! \bm{G} \bm{E}^T\! \dot{\bm{\theta}}  \!\!-\!\! \bm{B} \dot{\bm{\varphi}}_+ \!) \!\!-\!\! (\!\dot{\bm{q}}_{\rm fb} \!\!-\!\! \Delta \dot{\bm{q}}_{\rm fb}\!) \!\!\leq\!\!  ( \bm{1} \!\!-\!\! \dot{\bm{z}} ) \!M \label{eq-3-6-lpac:3}\\
            &  (\! \dot{\bm{z}} \!\!-\!\! \bm{1} \!) M  \!\! \leq \!\! 
            - \bm{B} \bm{1}  \!\!+\!\! (\! \bm{G} \bm{E}^T\! \dot{\bm{\theta}}  \!\!+\!\! \bm{B} \dot{\bm{\varphi}}_+ \!) \!\!-\!\! (\!\dot{\bm{q}}_{\rm tb} \!\!-\!\! \Delta \dot{\bm{q}}_{\rm tb}\!) 
            \!\!\leq\!\!  (\! \bm{1} \!\!-\!\! \dot{\bm{z}} ) \!M \label{eq-3-6-lpac:4}\\[-1mm]
            &  ( \dot{\bm{z}} \!-\! \bm{1} ) M  \! \leq \! - \bm{B} \bm{E}^T ( \dot{\bm{v}} - \bm{1}) \!-\! \Delta \dot{\bm{q}}_{\rm fb}  \!\leq\!   ( \bm{1} \!-\! \dot{\bm{z}} ) M  \label{eq-3-6-lpac:5}\\[-1mm]
            &   ( \dot{\bm{z}} \!-\! \bm{1} ) M  \! \leq \!~~  \bm{B} \bm{E}^T  ( \dot{\bm{v}} - \bm{1}) \!-\! \Delta \dot{\bm{q}}_{\rm tb} \!\leq\!  ( \bm{1} \!-\! \dot{\bm{z}} ) M   \label{eq-3-6-lpac:6}\\[-1mm]
            & \begin{aligned}
                & \bm{1}_{n_{\rm p}} \!\!\smallotimes\! \dot{\bm{\varphi}}_+ \!\!\leq\! \cos( \bm{\beta} \!\smallotimes\! \bm{\theta}_{\rm p} \!\!-\!\! \bm{1}_{n_{\rm p}} \!\!\smallotimes\! \bm{\theta}_{\rm max}  ) \\[-1mm]
                & -\! \sin( \bm{\beta} \!\smallotimes\! \bm{\theta}_{\!\rm p} \!\!-\!\! \bm{1}_{\!n_{\!\rm p}} \!\!\!\smallotimes\! \bm{\theta}_{\!\rm max}  ) \!\!\circ\!\! ( \bm{1}_{\!n_{\rm p}} \!\!\!\smallotimes\!\! \bm{E}^T \!\! \dot{\bm{\theta}} \!-\!\! \bm{\beta} \!\smallotimes\! \bm{\theta}_{\!\rm p} \!\!+\!\! \bm{1}_{\!n_{\rm p}} \!\!\smallotimes\! \bm{\theta}_{\!\rm max} ) 
            \end{aligned} \label{eq-3-6-lpac:7} \\[-1mm]
            & 0 \leq \dot{\bm{\varphi}}_+ \leq 1  \label{eq-3-6-lpac:8} 
        \end{align}  
\end{subequations}
where (\ref{eq-3-6-lpac:1}) and (\ref{eq-3-6-lpac:2}) model active power flow; (\ref{eq-3-6-lpac:3}) and (\ref{eq-3-6-lpac:4}) ((\ref{eq-3-6-lpac:5}) and (\ref{eq-3-6-lpac:6})) model the components of reactive power flow which are independent of (coupled with) changes in voltage magnitudes, i.e., $\bm{q}_{\rm fb} \!\!-\! \Delta \bm{q}_{\rm fb}$ and $\bm{q}_{\rm tb} \!\!-\! \Delta \bm{q}_{\rm tb}$ ($\Delta \bm{q}_{\rm fb}\!$ and $\Delta \bm{q}_{\rm tb}$); (\ref{eq-3-6-lpac:7}) defines the polyhedral outer approximation of the cosine function; (\ref{eq-3-6-lpac:8}) enforces limits on the cosine approximation; $\bm{\theta}_{\rm p} \!\!=\!\! 2 \bm{\theta}_{\rm max}/(1 \!+\! n_{\rm p})$ is the vector of segment length; and $\bm{\beta} \!\!=\!\! [1,2, \cdots\!, n_{\rm p}]^T$. 
In \cite{4-1163}, two cosine approximation variables to respectively approximate the associated entries in $\cos(\!\bm{E}^T \!\bm{\theta})$ and $\cos(\!-\! \bm{E}^T \!\bm{\theta})$, i.e., those in $\bm{\varphi}_{\!+}$  and $\bm{\varphi}_{\!-}$, are introduced for each branch; $\bm{\varphi}_+$ in (\ref{eq-3-6-lpac:4}) is replaced by $\bm{\varphi}_{\!-}$; and polyhedral outer approximation for $\bm{\varphi}_{\!-}$ analogous to (\ref{eq-3-6-lpac:7}) is contained. Ref. \cite{4-1129} uses a new approximation to reduce variables in $\bm{\varphi}_{\!+}$ and $\bm{\varphi}_{\!-}$ by half, resulting in substantial reductions of solution time. However, provided that the polyhedral outer approximation is symmetric about the $y$-axis, $\bm{\varphi}_-$ in \cite{4-1163} can be eliminated, as given by (\ref{eq-3-6-lpac}).

For (\ref{eq-3-6-ac-power-1:2}), we approximate $\dot{\bm{v}}^2$ by its first-order Taylor series at $\bm{1}$ as $\dot{\bm{v}}^2 \approx 2 \bm{v} - \bm{1} $. Substituting this approximation into (\ref{eq-3-6-ac-power-1:2}) and eliminating bilinear terms by the Big-M method yield
\begin{subequations}\label{eq-3-6-power-balance}  
    \begin{align}
        & [ \dot{\bm{v}} \bm{1}_{|\mathcal{E}|}^T \!-\! M (\bm{J}_{|\mathcal{V}| \!\times\! |\mathcal{E}|} \!-\! \bm{1}_{|\mathcal{V}|} \dot{\bm{z}}^T) ] \!\circ\! \tilde{\bm{E}} \!\leq\! \dot{\bm{U}} \!\leq\! M \bm{1}_{|\mathcal{V}|} \dot{\bm{z}}^T \!\circ\! \tilde{\bm{E}} \label{eq-3-6-power-balance:1}  \\
        & \!-\! M \bm{1}_{|\mathcal{V}|} \dot{\bm{z}}^T \!\!\!\circ\! \tilde{\bm{E}} \!\leq\! \dot{\bm{U}} \!\leq\! [\dot{\bm{v}} \bm{1}_{|\mathcal{E}|}^T \!+\! M (\!\bm{J}_{|\mathcal{V}| \!\times\! |\mathcal{E}|} \!-\! \bm{1}_{|\mathcal{V}|} \dot{\bm{z}}^T)] \!\circ\! \tilde{\bm{E}} \label{eq-3-6-power-balance:2}  \\
        & \bm{\Lambda} = 
        \begin{bmatrix}
            2 \dot{\bm{U}} \underline{\bm{G}}_{\rm b} \bm{1}_{|\mathcal{E}|}    \\
            - 2 \dot{\bm{U}} \underline{\bm{B}}_{\rm b} \bm{1}_{|\mathcal{E}|}   
        \end{bmatrix} 
        \!\!-\!\!
        \begin{bmatrix}
            \tilde{\bm{E}} \underline{\bm{G}}_{\rm b} \dot{\bm{z}} \\
            \!-\! \tilde{\bm{E}} \underline{\bm{B}}_{\rm b} \dot{\bm{z}}
        \end{bmatrix}  
        \!\!+\!\! 
        \begin{bmatrix}
            \underline{\bm{G}} \bm{1}   \!\!\\
            \!-\! \underline{\bm{B}} \bm{1} \!\!
        \end{bmatrix}
        \!\!\circ\!\! 
        \begin{bmatrix}
            2 \bm{v} - \bm{1} \\
            2 \bm{v} - \bm{1} 
        \end{bmatrix} 
    \end{align}
\end{subequations}
where $\dot{\bm{U}} \!\in\! \mathbb{R}^{|\mathcal{V}| \times |\mathcal{E}|}$ is an auxiliary variable matrix, and $\bm{\Lambda}$ represents the left-hand side (LHS) of (\ref{eq-3-6-ac-power-1:2}).

\subsubsection{Operational Constraints} They are formulated as 
\begin{subequations}\label{eq-3-6-operation-1}
\begin{align}
        & \bm{p}_{\rm c}^{\rm min} \leq \dot{\bm{p}}_{\rm c} \leq \bm{p}_{\rm c}^{\rm max}, \bm{q}_{\rm c}^{\rm min} \leq \dot{\bm{q}}_{\rm c} \leq \bm{q}_{\rm c}^{\rm max}  \label{eq-3-6-operation-1:1} \\[-1mm]
        &  \bm{0} \!\leq\!  \dot{\bm{p}}_{\rm v} \!\leq\!  \dot{\bm{p}}_{\rm v}^{\rm max}  \label{eq-3-6-operation-1:2}\\[-1mm]
        & \dot{\bm{q}}_{\rm v} \leq \dot{\bm{p}}_{\rm v} \circ \tan(\arccos ( \bm{\phi}_{\rm v}^{\rm min} ) )  \label{eq-3-6-operation-1:3}\\[-1mm]
        & \bm{v}_{\rm min}  \leq \dot{\bm{v}} \leq \bm{v}_{\rm max}  \label{eq-3-6-operation-1:4}\\[-1mm]
        & - \bm{\theta}_{\rm max} \!-\! (\bm{1} \!-\! \dot{\bm{z}}) M \leq \bm{E}^T \dot{\bm{\theta}} \leq  \bm{\theta}_{\rm max} \!+\! (\bm{1} \!-\! \dot{\bm{z}} ) M   \label{eq-3-6-operation-1:5}\\[-1mm]
        & \dot{\bm{p}}_{\rm v}^{ 2}  \!+\! \dot{\bm{q}}_{\rm v}^{ 2} \!\leq\! (\bm{s}_{\rm v}^{\rm max})^{ 2} \label{eq-3-6-operation-1:2-s}\\[-1mm]
        & \dot{\bm{p}}_{\rm fb}^{2} + \dot{\bm{q}}_{\rm fb}^{2}  \leq  \dot{\bm{z}} \circ  (\bm{s}_{\rm b}^{\rm max})^{2}, 
        \dot{\bm{p}}_{\rm tb}^{2} + \dot{\bm{q}}_{\rm tb}^{2}  \leq  \dot{\bm{z}} \circ  (\bm{s}_{\rm b}^{\rm max})^{2} 
        \label{eq-3-6-operation-1:6}
\end{align}
\end{subequations}
where (\ref{eq-3-6-operation-1:1}), (\ref{eq-3-6-operation-1:2}) and (\ref{eq-3-6-operation-1:2-s}) are output power constraints of generators, and VRE curtailment is allowed in (\ref{eq-3-6-operation-1:2}) considering potential VRE oversupply; 
(\ref{eq-3-6-operation-1:3}) limits the power factors of VRE-based generators to be greater than their minimum allowed values; 
(\ref{eq-3-6-operation-1:4}), (\ref{eq-3-6-operation-1:5}) and (\ref{eq-3-6-operation-1:6}) bound bus voltage magnitude, branch phase angle difference and branch power, respectively.

For constraints (\ref{eq-3-6-operation-1:6}), we use regular octagons as the outer approximation for the feasible regions \cite{4-1228, 4-1129}. Then, (\ref{eq-3-6-operation-1:6}) is approximated as 
\begin{subequations}\label{eq-3-6-first-stage-opera-2}
        \begin{align}
                &\!\!\!\!\! 
                - \bm{J}_2 \smallotimes (\dot{\bm{z}} \circ \bm{s}_{\rm b}^{\rm max})
                \!\!\leq\!\!
                \begin{bmatrix}
                        \dot{\bm{p}}_{\rm fb} &\!\!\!\! \dot{\bm{q}}_{\rm fb} \\
                        \dot{\bm{p}}_{\rm tb} &\!\!\!\! \dot{\bm{q}}_{\rm tb}
                \end{bmatrix} 
                \!\!\leq\!\! \bm{J}_2 \smallotimes (\dot{\bm{z}} \circ \bm{s}_{\rm b}^{\rm max}) \\
                &\!\!\!\!\! 
                -\!\!  \bm{J}_2 \!\smallotimes\! (\! \dot{\bm{z}} \!\circ\! \bm{s}_{\rm b}^{\rm max} \!)
                \!\!\leq\!\! 
                \frac{1}{\sqrt{2}} \!\! \begin{bmatrix}
                        \!\dot{\bm{p}}_{\rm fb} \!\!\!+\!\! \dot{\bm{q}}_{\rm fb} &\!\!\! \dot{\bm{p}}_{\rm fb} \!\!\!-\!\! \dot{\bm{q}}_{\rm fb}\!\! \\
                        \!\dot{\bm{p}}_{\rm tb} \!\!\!+\!\! \dot{\bm{q}}_{\rm tb} &\!\!\! \dot{\bm{p}}_{\rm tb} \!\!\!-\!\! \dot{\bm{q}}_{\rm tb}\!\!
                \end{bmatrix}
                \!\!\leq\!\!
                 \bm{J}_2 \!\smallotimes\! (\! \dot{\bm{z}} \!\circ\! \bm{s}_{\rm b}^{\rm max} \!)
        \end{align}
\end{subequations} 

Constraints (\ref{eq-3-6-operation-1:2-s}) can also be approximated similarly to (\ref{eq-3-6-operation-1:6}), but observing that $\bm{p}_{\rm v}$ and $\bm{q}_{\rm v}$ are also constrained by (\ref{eq-3-6-operation-1:2}) and (\ref{eq-3-6-operation-1:3}), (\ref{eq-3-6-operation-1:2-s}) only needs to be approximated at the part intersecting with the regions defined by (\ref{eq-3-6-operation-1:2}) and (\ref{eq-3-6-operation-1:3}). Specifically, (\ref{eq-3-6-operation-1:2-s}) is approximated as 
\begin{equation}\label{eq-3-6-first-stage-opera-3}
        \begin{aligned}
                &  \bm{1}_{n_{\rm s}} \smallotimes \dot{\bm{p}_{\rm v}} + \tan( \bm{\gamma} \smallotimes \bm{\theta}_{\rm s}  - \frac{\pi}{2} \!\cdot\! \bm{1}_{n_{\rm s} \cdot |\mathcal{U}_{\rm v}| }  ) \circ (\bm{1}_{n_{\rm s}} \smallotimes \dot{\bm{q}}_{\rm v} ) \\[-2mm] 
                & \leq \sec(   \bm{\gamma} \smallotimes \bm{\theta}_{\rm s}  - \frac{\pi}{2} \!\cdot\! \bm{1}_{n_{\rm s} \cdot |\mathcal{U}_{\rm v}| } ) \circ (\bm{1}_{n_{\rm s}} \smallotimes \bm{s}_{\rm v}^{\rm max} ) 
        \end{aligned}
\end{equation} 
where $\bm{\theta}_{\rm v}^{\rm max} \!=\!  \arccos ( \bm{\phi}_{\rm v}^{\rm min} \!)$, $\bm{\theta}_{\rm s} \!=\! (\frac{\pi}{2} \!+\! \bm{\theta}_{\rm v}^{\rm max})/(n_{\rm s} \!-\! 1)$, and $\bm{\gamma} \!=\! [0~ 1~ \cdots n_{\rm s} \!-\! 1]^T$.

\subsubsection{Topological Constraints}
Branches, e.g., transformers, do not participate in OTS as they stay switched on, giving
\begin{equation}\label{eq-3-6-remain-on}
        \bm{E}_{\rm on} \dot{\bm{z}} = \bm{1} 
\end{equation} 
where $\bm{E}_{\rm on}$ is the incidence matrix between branches remaining switched on and $\mathcal{E}$. In addition, network topology needs to be connected under the normal operation state, which can be ensured by the electrical flow-based constraints as follows:
\begin{equation}\label{eq-3-6-connected-1} 
  \begin{aligned}
    \!\!\! M \!(\!\dot{\bm{z}} \!-\! \bm{1} ) \!\!\leq\!\! \bm{E}^T\! {\bm{\vartheta}} \!\!-\!\! {\bm{\rho}}  \!\leq\! \! M\! (\!\bm{1} \!-\!\! \dot{\bm{z}} ), - M \dot{\bm{z}} \!\leq\! {\bm{\rho}} \!\leq\! M \dot{\bm{z}}, \bm{E} {\bm{\rho}} \!=\! \dot{\bm{\alpha}}
  \end{aligned}
\end{equation}
where $\bm{\rho} \!\!\in\!\! \mathbb{R}^{ |\mathcal{E}| }$ and $\bm{\vartheta} \!\!\in\!\! \mathbb{R}^{|\mathcal{V}|} $ are vectors of auxiliary variables, and $\dot{\bm{\alpha}}$ is a $|\mathcal{V}|$-dimensional constant uniquely-balanced vector. 
Constraints (\ref{eq-3-6-connected-1}) ensure NC of the power grid by leveraging the equivalence between it and feasibility of the vertex potential equation of a  specially designed electrical flow network. 
For more details of (\ref{eq-3-6-connected-1}) we refer the reader to \cite{4-995-ea}.

The OTS model is symmetric if multiple identical lines exist between any two buses \cite{4-1154}. Adding static symmetry-breaking inequalities can cut the redundant symmetric solutions thus reducing solution time potentially \cite{4-1239}. These inequalities are$\!\!$
\begin{equation}\label{eq-3-6-symmetry}
       \bm{K}_{(|s|-1) \times |s|} \bm{E}^{{\rm sym}}_s  \dot{\bm{z}} \leq \bm{0}   ~~ \forall s \in \mathcal{S}
\end{equation}
where $\bm{E}^{{\rm sym}}_s$ is the incidence matrix between $s$ and $\mathcal{E}$. It is noted that for the identical lines, every property involved in the OTS model needs to be identical. 


\subsection{Second-Stage Formulation}

\subsubsection{Objective} The second-stage objective is to minimize the regulation cost of active power outputs, i.e., 
\begin{equation}\label{eq-3-6-second-stage-obj-linear}
        \min~
        \ddot{\bm{w}}_{\rm c_+}^T \ddot{\bm{p}}_{\rm c_+} 
        \!+\! \ddot{\bm{w}}_{\rm c_-}^T \ddot{\bm{p}}_{\rm c_-}   
        + \ddot{\bm{w}}_{\rm v_+}^T \ddot{\bm{p}}_{\rm v_+} 
        \!+\! \ddot{\bm{w}}_{\rm v_-}^T \ddot{\bm{p}}_{\rm v_-}  
\end{equation} 
Generally, we have $\bm{\omega}_{\rm c_+} \!\!\gg\!\! \bm{\omega}_{\rm v_+}$ and $\bm{\omega}_{\rm c_-} \!\!\ll\!\! \bm{\omega}_{\rm v_-}$ to promote utilization of VRE and also penalize curtailment of VRE.

\subsubsection{Constraints} 
The second-stage constraints include all first-stage power flow and operational constraints with first-stage variables and parameters replaced by second-stage ones except for $\dot{\bm{z}}$ and $\dot{\bm{p}}_{\rm d}$, and additional constraints as follows:
\begin{subequations}\label{eq-3-6-second-stage-cst}
    \begin{align}
        & \bm{0} \leq \ddot{\bm{p}}_{\rm c_+} \leq  \Delta \ddot{t} \!\cdot\! {\bm{r}}_{\rm c_+}, \bm{0}  \leq  \ddot{\bm{p}}_{\rm c_-}  \leq  \Delta \ddot{t} \!\cdot\! {\bm{r}}_{\rm c_-} 
        \label{eq-3-6-second-stage-cst:1}\\
        & \bm{0} \!\!\leq\!\! \ddot{\bm{p}}_{\rm v_{\!+}} \!\!\leq\!\! \Delta \ddot{t} \!\cdot\! {\bm{r}}_{\!\rm v_{\!+}}, \bm{0}  \!\!\leq\!\! \min(\dot{\bm{p}}_{\rm v}, \ddot{\bm{p}}_{\rm v}^{\rm max})  \!+\! \ddot{\bm{p}}_{\rm v_{\!+}} \!\!-\! \ddot{\bm{p}}_{\rm v}  \!\!\leq\!\! \Delta \ddot{t} \!\cdot\! {\bm{r}}_{\rm v_{\!-}} 
        \label{eq-3-6-second-stage-cst:2}\\
        & \ddot{\bm{p}}_{\rm v}^{\rm max} = \dot{\bm{p}}_{\rm v}^{\rm max} + \bm{\varepsilon}
        \label{eq-3-6-second-stage-cst:3} \\
        &  \ddot{\bm{p}}_{\rm c} =  \dot{\bm{p}}_{\rm c} + \ddot{\bm{p}}_{\rm c_+} - \ddot{\bm{p}}_{\rm c_-}, \ddot{\bm{p}}_{\rm v} =  \dot{\bm{p}}_{\rm v} + \ddot{\bm{p}}_{\rm v_+} - \ddot{\bm{p}}_{\rm v_-}  \label{eq-3-6-second-stage-cst:4}\\
        & \ddot{\bm{v}}_{\rm c} = \dot{\bm{v}}_{\rm c} + \ddot{\bm{v}}_{\rm c_{\Delta}}, \ddot{\bm{q}}_{\rm v} = \dot{\bm{q}}_{\rm v} + \ddot{\bm{q}}_{\rm v_{\Delta}} \label{eq-3-6-second-stage-cst:5}
    \end{align}
\end{subequations}
where (\ref{eq-3-6-second-stage-cst:1}) and (\ref{eq-3-6-second-stage-cst:2}) are ramping constraints of conventional and VRE-based generators respectively, imposed on the active power output regulation from the first-stage scheduled value to the second-stage one by corrective controls; 
(\ref{eq-3-6-second-stage-cst:3}) gives the true values of the maximal active power outputs of VRE-based generators, $\ddot{\bm{p}}_{\rm v}^{\rm max} $, by its forecast $\dot{\bm{p}}_{\rm v}^{\rm max}$ and forecast error $\bm{\varepsilon}$;
(\ref{eq-3-6-second-stage-cst:4}) and (\ref{eq-3-6-second-stage-cst:5}) are equalities for variables changed directly by the second-stage corrective controls. Note that VRE-based generators are assumed to be equipped with control loops to limits ramp rates to fulfill the grid code \cite{4-1328, 4-1330}, giving ramping constraints (\ref{eq-3-6-second-stage-cst:2}). Moreover, the second part of (\ref{eq-3-6-second-stage-cst:2}) can be easily linearized via introducing auxiliary continuous variables to represent the difference between $\dot{\bm{p}}_{\rm v}$ and $\ddot{\bm{p}}_{\rm v}^{\rm max}$ and penalizing the sum of the introduced variables in (\ref{eq-3-6-second-stage-obj-linear}).

\subsection{Third-Stage Formulation}


\subsubsection{Objective} The third-stage objective is to minimize the corrective control cost for responding to contingencies, i.e., 
\begin{equation}\label{eq-3-6-third-stage-obj-linear}
     \!\!\!\!\!\!\!\! 
     \min~  
        \dddot{\bm{w}}_{\!\rm c_{\!+}}^T \! \dddot{\bm{p}}_{\!\rm c_{\!+}} 
        \!\!\!+\!\! \dddot{\bm{w}}_{\!\rm c_{\!-}}^T \! \dddot{\bm{p}}_{\!\rm c_{\!-}}  
        \!\!\!+\!\! \dddot{\bm{w}}_{\!\rm v_{\!+}}^T \! \dddot{\bm{p}}_{\!\rm v_{\!+}} 
        \!\!\!+\!\! \dddot{\bm{w}}_{\!\rm v_{\!-}}^T \! \dddot{\bm{p}}_{\!\rm v_{\!-}}
        \!\!\!+\!\! \bm{w}_{\!\rm d}^T \! \bm{p}_{\rm d_{\!\Delta}}
        \!\!\!+\!\! \bm{w}_{\!\rm s}^T \!(\!\bm{z}_{\!+} \!\!+\!\! \bm{z}_{\!-}\!) \!\!\!\!\!\! 
\end{equation}
where the first four terms give the regulation cost of generators, and the last two are the load shedding cost and line switching cost, respectively.

\subsubsection{Constraints} In addition to the first-stage power flow constraints, (\ref{eq-3-6-operation-1:4}), (\ref{eq-3-6-operation-1:5}), (\ref{eq-3-6-first-stage-opera-2}), and (\ref{eq-3-6-first-stage-opera-3}), all with the first-stage variables replaced by the third-stage ones, the third-stage constraints also include 
\begin{subequations}\label{eq-3-6-3-stage-opera}
    \begin{align}
       & 
       \bm{0} \!\leq\! \dddot{\bm{p}}_{\rm c_+} \!\leq\!  \bm{o}_{\rm c} \!\circ\! (\Delta\! \dddot{t} \cdot {\bm{r}}_{\rm c_+}), 
       \bm{0}  \!\leq\!  \dddot{\bm{p}}_{\rm c_-}  \!\!\leq\!   \bm{o}_{\rm c} \circ (\Delta\! \dddot{t} \cdot{\bm{r}}_{\rm c_-}) 
       \label{eq-3-6-3-stage-opera:1}\\[-1mm]
       & 
       \bm{0} \!\leq\! \dddot{\bm{p}}_{\rm v_+} \!\leq\!   \bm{o}_{\rm v} \!\circ\! (\Delta\! \dddot{t} \cdot {\bm{r}}_{\rm v_+}), 
       \bm{0}  \!\leq\!  \dddot{\bm{p}}_{\rm v_-}  \!\!\leq\!  \bm{o}_{\rm v} \!\circ\! (\Delta\! \dddot{t} \cdot {\bm{r}}_{\rm v_-}) 
       \label{eq-3-6-3-stage-opera:2}\\[-1mm]
       & \bm{o}_{\rm c} \!\circ\! \bm{p}_{\rm c}^{\rm min} \!\leq\! \dddot{\bm{p}}_{\rm c} \leq \bm{o}_{\rm c} \!\circ\! \bm{p}_{\rm c}^{\rm max}, \bm{o}_{\rm c} \!\circ\! \bm{q}_{\rm c}^{\rm min} \!\leq\! \dddot{\bm{q}}_{\rm c} \!\leq\! \bm{o}_{\rm c} \!\circ\! \bm{q}_{\rm c}^{\rm max}  
       \label{eq-3-6-3-stage-opera:3}\\[-1mm]
       &  \bm{0} \!\leq\!  \dddot{\bm{p}}_{\rm v} \!\leq\!  \bm{o}_{\rm v} \!\circ\! \ddot{\bm{p}}_{\rm v}^{\rm max} 
       \label{eq-3-6-3-stage-opera:4}\\[-1mm] 
       & - \bm{o}_{\rm v} \circ \bm{s}_{\rm v}^{\rm max} \leq \dddot{\bm{q}}_{\rm v} \leq \bm{o}_{\rm v} \circ  \dddot{\bm{p}}_{\rm v} \circ \tan(\arccos ( \bm{\phi}_{\rm v}^{\rm min} ) ) 
       \label{eq-3-6-3-stage-opera:5} \\[-1mm]
       &  \dddot{\bm{p}}_{\rm c} \!=\!   \bm{o}_{\rm c} \!\circ\! \ddot{\bm{p}}_{\rm c} \!+\! \dddot{\bm{p}}_{\rm c_+} \!-\! \dddot{\bm{p}}_{\rm c_-}, \dddot{\bm{p}}_{\rm v} \!=\! \bm{o}_{\rm v} \!\circ\! \ddot{\bm{p}}_{\rm v} \!+\! \dddot{\bm{p}}_{\rm v_+} \!-\! \dddot{\bm{p}}_{\rm v_-}  \label{eq-3-6-3-stage-opera:6}\\[-1mm]
       & \dddot{\bm{v}}_{\rm c} =   \ddot{\bm{v}}_{\rm c} + \dddot{\bm{v}}_{\rm c_{\Delta}}, \dddot{\bm{q}}_{\rm v} = \ddot{\bm{q}}_{\rm v} + \dddot{\bm{q}}_{\rm v_{\Delta}},
       \dddot{\bm{p}}_{\rm d} = \dot{\bm{p}}_{\rm d} - \bm{p}_{\rm d_{\Delta}}
       \label{eq-3-6-3-stage-opera:7} \\[-1mm]
       &  \bm{0} \leq \bm{p}_{\rm d_{\Delta}} \leq  \bm{p}_{\rm d_{\Delta}}^{\rm max}, 0 \leq \dddot{\bm{q}}_{\rm d} \leq \dot{\bm{q}}_{\rm d}, \dddot{\bm{p}}_{\rm d} \circ \dot{\bm{q}}_{\rm d} = \dddot{\bm{q}}_{\rm d} \circ \dot{\bm{p}}_{\rm d}  
       \label{eq-3-6-3-stage-opera:8} 
    \end{align}
\end{subequations}
where (\ref{eq-3-6-3-stage-opera:1}) to (\ref{eq-3-6-3-stage-opera:5}) are ramping and power output constraints for generators similarly to those for the second stage while with parameterization of generator contingencies such that power outputs and regulations of fault generators are forced to 0; 
(\ref{eq-3-6-3-stage-opera:6}) and (\ref{eq-3-6-3-stage-opera:7}) are equalities for variables changed directly by the third-stage corrective controls; 
and (\ref{eq-3-6-3-stage-opera:8}) limits load shedding amounts and keeps constant load power factors during load shedding.

\subsubsection{Topological Constraints} 
\begin{figure}[h]
	\centering
	\includegraphics[scale=0.9]{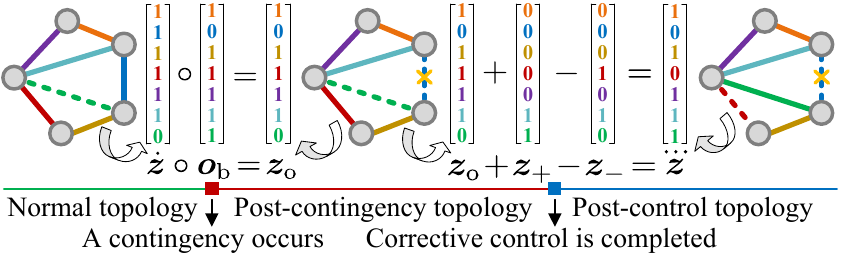}  
	\caption{Different topologies and their relationship.}
	\label{fig-3-6-4}
\end{figure}

Two topologies can appear at the third stage: the topology after contingencies but before corrective controls (post-contingency topology for short), denoted by $\bm{z}_{\rm o}$, and that after corrective controls (post-control topology for short), denoted by $\dddot{\bm{z}}$. 
Fig. \ref{fig-3-6-4} illustrates these two topologies and the normal topology at the first stage, and the relationship between them using a 5-bus network. In this network, the green line is switched off at the first stage, and then the blue line is opened by a contingency, followed by the green and blue lines being switched on and off by the corrective control respectively.

The post-contingency topology is formed by imposing the branch contingency parameterized by $\bm{o}_{\rm b}$ on the normal topology $\dot{\bm{z}}$, which gives
\begin{equation}\label{eq-3-6-third-stage-con-z-0}
    \bm{z}_{\rm o} = \dot{\bm{z}} \!\circ\! \bm{o}_{\rm b}
\end{equation}

The post-control topology $\dddot{\bm{z}}$ is obtained by imposing the corrective control parameterized by $\bm{z}_{+}$ and $\bm{z}_{-}$ on the post-contingency topology $\bm{z}_{\rm o}$. Thus we have
\begin{equation}\label{eq-3-6-third-stage-con-z-a}
    \dddot{\bm{z}} = \bm{z}_{\rm o} + \bm{z}_+ - \bm{z}_-
\end{equation}

Variables $\bm{z}_+$ and $\bm{z}_{-}$ that parameterize the corrective switching actions are binary variables. However, to reduce the computational complexity, instead of explicitly defining $\bm{z}_+$ and $\bm{z}_{-}$ as binary variables, we define $\bm{z}_+$ and $\bm{z}_{-}$ as continuous variable and introduce the following constraints:
\begin{equation}\label{eq-3-6-third-stage-con-z-1}
    \bm{0} \leq \bm{z}_+ \leq \bm{1}, \bm{0} \leq \bm{z}_- \leq \bm{1}, \bm{z}_+ + \bm{z}_- \leq \bm{1} 
\end{equation}
where the first two inequalities bound $\bm{z}_+$ and $\bm{z}_-$, respectively; and the last one is obtained with the fact that one branch should not be both switched on and off by the corrective control. Since $\dddot{\bm{z}}$ and $\bm{z}_{\rm o}$ are binary variables, equality (\ref{eq-3-6-third-stage-con-z-a}) implies that $\bm{z}_+ - \bm{z}_-$ are also binary. This together with constraints (\ref{eq-3-6-third-stage-con-z-1}) and the fact that $\bm{w}_{\rm s}^T (\bm{z}_{+} \!+\! \bm{z}_{-})$ is minimized in the objective function (\ref{eq-3-6-third-stage-obj-linear}), ensures that $\bm{z}_+$ and $\bm{z}_{-}$ are binary variables.

Branches with failures should be open in the post-control topology, and branches that do not participate in the corrective controls execute no switching actions. Thus we have
\begin{equation}\label{eq-3-6-third-stage-con-z-2}
   \dddot{\bm{z}} \leq \bm{o}_{\rm b}, \bm{E}_{\rm nc} \bm{z}_+ = \bm{0}, \bm{E}_{\rm nc} \bm{z}_- = \bm{0}
\end{equation}
where $\bm{E}_{\rm nc}$ is the incidence matrix between branches that do not participate in the corrective controls and $\mathcal{E}$. 


Next, we consider NC at the third stage, based on the network connectedness criteria for security-constrained OTS and associated formulation developed by us in \cite{4-1355}. 
To begin with, we introduce a parameterized region $\mathcal{Z}( \phi,\! \dddot{\bm{\alpha}},\! \tilde{\bm{\alpha}} )\!=\!$
\begin{equation}\label{eq-3-6-connected-uniform-defination}
     \!\!\!\!\!\! 
     \left\{ \!\!
      \bm{z} \!\!\in\!\! \mathbb{B}^{|\mathcal{E}|} \!\left\vert \!\!
    \begin{aligned} 
        & M (\bm{z}  \!\!-\!\! \bm{1} ) \!-\!  \phi M \bm{1} \!\!\leq\!\! \bm{E}^T \!\tilde{\bm{\vartheta}} \!-\!\! \tilde{\bm{\rho}}  \!\leq\! M (\bm{1} \!\!-\!\! \bm{z} ) \!+\!   \phi M \bm{1}  \\  
        & \!-\! M \bm{z} \!-\! \phi M \bm{1} \!\leq\! \tilde{\bm{\rho}} \!\leq\! M \bm{z} \!+\! \phi M \!\bm{1}  \\
        & -\! \phi M \bm{1} \!\!\leq\!\! \bm{E} \tilde{\bm{\rho}} \!-\!\! (\dddot{\bm{\alpha}} \!+\! \tilde{\bm{\alpha}}) \!\!\leq\!\! \phi M \!\bm{1},\! 
        \tilde{\bm{\vartheta}} \!\!\in\!\! \mathbb{R}^{|\mathcal{V}|}\!, \tilde{\bm{\rho}} \!\!\in\!\! \mathbb{R}^{|\mathcal{V}|} 
    \end{aligned} 
     \right. \!
    \!\right\} \!\!\!\!
\end{equation}
with $\phi \in \mathbb{R}$, $\tilde{\bm{\alpha}} \in \mathbb{R}^{|\mathcal{V}|}$ and $\dddot{\bm{\alpha}}$ is a $|\mathcal{V}|$-dimensional constant uniquely-balanced vector \cite{4-1355}. When $\phi\!=\!0$ and $\tilde{\bm{\alpha}} \!\!=\!\! \bm{0}$, $\mathcal{Z}( \phi,\! \dddot{\bm{\alpha}},\! \tilde{\bm{\alpha}} )$ is the region of $\bm{z}$ corresponding to all connected topologies, which is imposed by constraints analogous to (\ref{eq-3-6-connected-1}). When $\phi\!\geq\!1$ and $\tilde{\bm{\alpha}} \!\!=\!\! \bm{0}$, all constraints to impose $\bm{z}$ are invalid. 
Furthermore, let $(\bm{\alpha}_+^*, \bm{\alpha}_-^*)$ be the optimal value of $(\bm{\alpha}_+, \bm{\alpha}_-)$ of the following linear programming model:
\begin{subequations}\label{eq-3-6-connected-3}
    \begin{align}
        \min\nolimits_{\bm{\alpha}_+ \in \mathbb{R}^{|\mathcal{V}|}\!, \bm{\alpha}_- \in \mathbb{R}^{|\mathcal{V}|} \!} & \bm{1}^T (\bm{\alpha}_+ + \bm{\alpha}_-) \\
    \text{s.t.}~ &  \bm{z}_{\rm o} \in \mathcal{Z}(0, \dddot{\bm{\alpha}}, \bm{\alpha}_+ \!-\! \bm{\alpha}_-  )  \label{eq-3-6-connected-3:1} \\
    & \bm{\alpha}_+ \geq \bm{0},\bm{\alpha}_- \geq \bm{0}  \label{eq-3-6-connected-3:3} 
    \end{align}
\end{subequations} 
Then, the following constraints are imposed on $\bm{z}_{\rm o}$ and $\dddot{\bm{z}}$:
\begin{subequations}\label{eq-3-6-topology-all}
    \begin{align}
        & 1 \!\!-\! \phi_1 (k_b \!+\! 1)  \!\leq\!  (|\mathcal{E}| \!-\! \bm{1}^T \!\bm{o}_{\rm b}) \!-\! k_{\rm b} \!\leq\!  (1 \!\!-\! \phi_1) (k_{\rm max} \!-\! k_{\rm b} ) \label{eq-3-6-topology-kb} \\
        & \bm{1}^T ( \bm{\alpha}_+^* + \bm{\alpha}_-^* ) \leq n_{\rm m} r + (1 - \phi_1) M \label{eq-3-6-connected-zo} \\
        & M \phi_3 \geq 2 \Vert \tilde{\bm{J}} \dddot{\bm{\alpha}} \Vert_{- 2 \infty}  - \bm{1}^T (\bm{\alpha}_+^* + \bm{\alpha}_-^*) \geq M (\phi_3 - 1) \label{eq-3-6-phi-3-4} \\
        & M \phi_4 \geq \bm{1}^T (\bm{\alpha}_+^* + \bm{\alpha}_-^*) - n_{\rm m} r \geq M (\phi_4 - 1) \label{eq-3-6-phi-3-4-1} \\
        & (\phi_2 - 1) k_{\rm max}  \leq  \bm{1}^T (\dot{\bm{z}}  - \bm{z}_{\rm o}) \leq ( 1 - \phi_2 ) k_{\rm max} \label{eq-3-6-topology-tau} \\
        & 1 - \phi_5  \leq  \bm{1}^T \bm{z}_-  \leq (1 - \phi_5)( |\mathcal{E}| - |\mathcal{V}| + 1 )  \label{eq-3-6-phi-5} \\
        & \dddot{\bm{z}} \in \mathcal{Z}([\bm{1}^T (\bm{\alpha}_+^* + \bm{\alpha}_-^*) + \phi_5], \dddot{\bm{\alpha}}, \bm{0} ) \label{eq-3-6-connected-5:1} \\
        & \dddot{\bm{z}} + (\bm{1} - \bm{o}_{\rm b}) \in \mathcal{Z}([\phi_3 + \phi_4], \dddot{\bm{\alpha}}, \bm{0} ) \label{eq-3-6-connected-5:2}
    \end{align}
\end{subequations}
where $k_{\max}$ is the maximal number of fault components; $k_{\rm b}$, $r$, $n_{\rm m}$ and $\tilde{\bm{J}}$ are all pre-defined parameters; $\phi_1$ to $\phi_5$ are all binary scalar variables. Variable $\phi_1$ indicates if the number of fault branches of the contingency is not more than $k_{\rm b}$, $\phi_2$ indicates if $\bm{z}_{\rm o}$ equals to $\dot{\bm{z}}$, $\phi_3$ and $\phi_4$ together indicate if $\bm{z}_{\rm o}$ is unconnected with only inevitable disconnection, and $\phi_5$ indicates if any lines are switched off in corrective controls. 
Constraints (\ref{eq-3-6-topology-all}) together guarantee that the post-contingency topology is connected for all contingencies where the number of fault branches is not more than $k_{\rm b}$, ignoring inevitable network disconnection; and corrective line switching does not further disconnect the network when the post-contingency topology is connected or disconnected with only inevitable disconnection.




\subsection{Contingency Uncertainties}

Denote by $\Omega$ the support of $\bm{o}$, $\mathcal{F}$ its $\xi$-field, and $\mathcal{M}(\Omega, \mathcal{F})$ the set of all PDs defined on $(\Omega, \mathcal{F})$. Then $\Omega$ is given by 
\begin{equation} 
         \Omega = \{ \bm{o} \in \mathbb{B}^N: N - \bm{1}^T \bm{o} \leq k_{\max}  \} 
\end{equation} 
An accurate PD of components failures is generally unattainable since the failure probabilities are close to 0 and most failures have not happened in historical data \cite{4-1214}. Thus, the following moment-based ambiguity set is adopted \cite{4-1214}:
\begin{equation}\label{eq-3-6-ambiguity-set-1}
    \mathcal{O} = \{ O \in \mathcal{M}(\Omega, \mathcal{F}):  {\bm{o}}_{\rm min}  \leq   \mathbb{E}_{O}( \bm{1} - \bm{o} )   \!\leq\! {\bm{o}}_{\rm max}    \}
\end{equation}
where ${\bm{o}}_{\rm min}$ and ${\bm{o}}_{\rm max}$ are vectors of lower and upper bounds of probabilities of components failures, respectively. Interval estimation can provide the confidence intervals of probabilities of components failures for a designated confidence level. In practice, ${\bm{o}}_{\rm min}$ and ${\bm{o}}_{\rm max}$ can be set as the lower and upper bounds of confidence intervals. Concisely, we rewrite (\ref{eq-3-6-ambiguity-set-1}) as
\begin{equation}\label{eq-3-6-ambiguity-set}
    \mathcal{O} = \{ O \in \mathcal{M}(\Omega, \mathcal{F}): \mathbb{E}_{O}( \bm{T} (\bm{1} - \bm{o}) )  \leq \tilde{\bm{o}}  \}
\end{equation}
where $\bm{T} \!=\! [ \bm{I}_{N} ~ -\!\bm{I}_{N} ]^T $ and $\tilde{\bm{o}} \!=\! [{\bm{o}}_{\max}^T ~ -\!{\bm{o}}_{\min}^T]^T$.


\subsection{Final Formulation}

The final TSDR-OTS model can be written in the following compact form:
\begin{subequations}\label{eq-3-6-final-stage-1}
    \begin{align}
        Q_1 := \min\nolimits_{\bm{x} \in \mathbb{R}^{n_1} \!\times \mathbb{B}^{m_1}} ~&
        \bm{c}_{\rm 1}^T \bm{x} +   \mathbb{E}_{{P}}[ Q_{\rm 2}( \bm{x}, \bm{\varepsilon} ) ] \label{eq-3-6-final-stage-1:1}\\
        \text{s.t.} ~& 
        \bm{A} \bm{x} \leq \bm{b}  \label{eq-3-6-final-stage-1:2} 
    \end{align}
\end{subequations}
with the second-stage recourse function
\begin{subequations}\label{eq-3-6-final-stage-2}
    \begin{align}
        \!\!\!\! Q_2( \bm{x}, \bm{\varepsilon} ) \!:=\! 
        \min_{\bm{x}_{\bm{\varepsilon}} \in \mathbb{R}^{n_2}} ~ \!\!& \bm{c}_{\rm 2}^T \bm{x}_{\bm{\varepsilon}} 
        \!+\! 
        \sup_{ O \in \mathcal{O} }   \mathbb{E}_{ O }   
            [ Q_3( \bm{x}_{\bm{\varepsilon}},\! \bm{x},\! \bm{o} ) ]  \label{eq-3-6-final-stage-2:1} \\
        \text{s.t.} ~ & 
        \bm{C} \bm{x}_{\bm{\varepsilon}} + \bm{D} \bm{x} \geq  \bm{d}( \bm{\varepsilon} )   \label{eq-3-6-final-stage-2:2} 
    \end{align}
\end{subequations}
and the third-stage recourse function
\begin{subequations}\label{eq-3-6-final-stage-3}
    \begin{align}
        \!\!\!\! Q_3(\! \bm{x},\! \bm{x}_{\bm{\varepsilon}},\! \bm{o} \!)   &  \!:=\! 
        \min\nolimits_{\bm{x}_{\bm{o}}^{\bm{\varepsilon}} \in \mathbb{R}^{n_3} \!\times \mathbb{B}^{m_3}} ~   \bm{c}_{\rm 3}^T \bm{x}_{\bm{o}}^{\bm{\varepsilon}}   \label{eq-3-6-final-stage-3:1}\\
        \text{s.t.} ~&  \bm{F} \bm{x}_{\bm{o}}^{\bm{\varepsilon}} \!+\! \bm{H} \bm{x}_{\bm{\varepsilon}} \!+\!  \bm{L} \bm{x} \!\geq\! \bm{f} \label{eq-3-6-final-stage-3:2} \\
        \!\!\!~ &  \bm{M} \bm{x}_{\bm{o}}^{\bm{\varepsilon}} + \bm{N} \bm{o} \geq \bm{g}  \label{eq-3-6-final-stage-3:3} \\
        \!\!\!~ & \bm{R} \bm{x}_{\bm{o}}^{\bm{\varepsilon}} \!+\! \bm{S} (\bm{o}) [\bm{x}^T ~\bm{x}_{\bm{\varepsilon}}^T ~(\bm{x}_{\bm{o}}^{\bm{\varepsilon}})^T ]^T   \!\geq\! \bm{h} \label{eq-3-6-final-stage-3:4}\\
        \!\!\!~ & [\bm{x}_{\bm{o}}^{\bm{\varepsilon}}]_{\rm c} \!\!\in\! \argmin_{ [\bm{x}_{\bm{o}}^{\bm{\varepsilon}}]_{\rm c} }  \{\! \bm{l}^T \![\bm{x}_{\bm{o}}^{\bm{\varepsilon}}]_{\rm c}\!\!:\! {  \bm{V} \![\bm{x}_{\bm{o}}^{\bm{\varepsilon}}]_{\rm z} \!\!+\!\! \bm{W} \![\bm{x}_{\bm{o}}^{\bm{\varepsilon}}]_{\rm c} \!\!\leq\!\! \bm{u}} \!\} \label{eq-3-6-final-stage-3-5}  \!\!\!
    \end{align}
\end{subequations} 
where 
$\bm{x}$, $\bm{x}_{\bm{\varepsilon}}$ and $\bm{x}_{\bm{o}}^{\bm{\varepsilon}}$ are optimization variables for different stages, including all control and state variables involved in each stage; 
$[\bm{x}_{\bm{o}}^{\bm{\varepsilon}}]_{\rm c}$ denotes the sub-vector of $\bm{x}_{\bm{\varepsilon}}$ corresponding to optimization variables of (\ref{eq-3-6-connected-3}), and $[\bm{x}_{\bm{o}}^{\bm{\varepsilon}}]_{\rm z}$ denotes the sub-vector of $\bm{x}_{\bm{o}}^{\bm{\varepsilon}}$ corresponding to $\bm{z}_{\rm o}$; 
$n_1$ (or $n_3$) and $m_1$ (or $m_3$) are dimensions of continuous variables and binary variables in $\bm{x}$ (or $\bm{x}_{\bm{o}}^{\bm{\varepsilon}}$), respectively, and $n_2$ the dimension of $\bm{x}_{\bm{\varepsilon}}$; 
$\bm{A}$, $\bm{C}$, $\bm{D}$, $\bm{F}$, $\bm{H}$, $\bm{L}$, $\bm{M}$, $\bm{N}$, $\bm{R}$, $\bm{S}$, $\bm{V}$ and $\bm{W}$ are proper matrices; $\bm{b}$, $\bm{f}$, $\bm{g}$, $\bm{h}$, $\bm{u}$, $\bm{l}$, and $\bm{c}_i$ are proper vectors; 
$\bm{d}( \bm{\varepsilon} )$ and $\bm{S}(\bm{o})$ denote a proper variable vector and matrix linearly related to $\bm{\varepsilon}$ and $\bm{o}$, respectively; 
$\bm{c}_{\rm 1}^T \bm{x}_{\rm 1}$, $\bm{c}_{\rm 2}^T \bm{x}_{\bm{\varepsilon}}$, and $\bm{c}_{\rm 3}^T \bm{x}_{\bm{o}}^{\bm{\varepsilon}}$ are equivalent to functions in (\ref{eq-3-6-first-stage-obj-linear}), (\ref{eq-3-6-second-stage-obj-linear}), and $\!$(\ref{eq-3-6-third-stage-obj-linear}), respectively; detailed formulations of the constraints are given in Table \ref{table-3-6-0}. 
It is noted that the expectation in (\ref{eq-3-6-final-stage-1}) can be tackled using the scenario-based approximation method as in \cite{4-1139}.

\begin{table}[h]
	\centering
    \caption{Detailed formulations of the constraints in (\ref{eq-3-6-final-stage-1})-(\ref{eq-3-6-final-stage-3})}
    \setlength{\tabcolsep}{0pt} 

    \setlength{\aboverulesep}{0pt}
    \setlength{\belowrulesep}{0pt}
    \setlength{\extrarowheight}{.1ex}
     \small{
    \begin{tabular*}{\hsize}{@{}p{1.5cm}@{}p{7.28cm} }\hline\hline 
               Constraint  & Detailed formulation    \\ \hline\hline
               (\ref{eq-3-6-final-stage-1:2}) 
               & 
               (\ref{eq-3-6-first-stage-obj-liearcons}), 
               (\ref{eq-3-6-lpac}), 
               (\ref{eq-3-6-power-balance}), 
               (\ref{eq-3-6-operation-1:1})-(\ref{eq-3-6-operation-1:5}), 
               (\ref{eq-3-6-first-stage-opera-2})-(\ref{eq-3-6-symmetry})   
               \\ \hline
               (\ref{eq-3-6-final-stage-2:2}) 
               & 
               \{(\ref{eq-3-6-lpac}), 
               (\ref{eq-3-6-power-balance}), 
               (\ref{eq-3-6-operation-1:1})-(\ref{eq-3-6-operation-1:5}), 
               (\ref{eq-3-6-first-stage-opera-2}), (\ref{eq-3-6-first-stage-opera-3})\} with all first-stage variables and parameters replaced by the second-stage ones except for $\dot{\bm{z}}$ and $\dot{\bm{p}}_{\rm d}$, (\ref{eq-3-6-second-stage-cst})   
               \\ \hline
               (\ref{eq-3-6-final-stage-3:2}) 
               & 
               \{(\ref{eq-3-6-lpac}), 
                (\ref{eq-3-6-power-balance}), 
                (\ref{eq-3-6-operation-1:4}), 
                (\ref{eq-3-6-operation-1:5}), 
                (\ref{eq-3-6-first-stage-opera-2}), 
                (\ref{eq-3-6-first-stage-opera-3})\} with all first-stage variables replaced by the third-stage ones, 
                the LHS in \{ (\ref{eq-3-6-3-stage-opera:1}), (\ref{eq-3-6-3-stage-opera:2}), (\ref{eq-3-6-3-stage-opera:4})\}, 
                (\ref{eq-3-6-3-stage-opera:7}), 
                (\ref{eq-3-6-3-stage-opera:8}), 
                (\ref{eq-3-6-third-stage-con-z-a}), 
                (\ref{eq-3-6-third-stage-con-z-1}), 
                the last two of (\ref{eq-3-6-third-stage-con-z-2}), 
                (\ref{eq-3-6-connected-zo}), 
                (\ref{eq-3-6-topology-tau}), 
                (\ref{eq-3-6-phi-3-4}), (\ref{eq-3-6-phi-3-4-1}), (\ref{eq-3-6-phi-5}), (\ref{eq-3-6-connected-5:1})  
                \\ \hline
                (\ref{eq-3-6-final-stage-3:3})
                & 
                The right-hand side (RHS) of (\ref{eq-3-6-3-stage-opera:1}) and (\ref{eq-3-6-3-stage-opera:2}), 
                the LHS of (\ref{eq-3-6-3-stage-opera:5}), 
                (\ref{eq-3-6-3-stage-opera:3}), 
                the first one in (\ref{eq-3-6-third-stage-con-z-2}), 
                (\ref{eq-3-6-topology-kb}), 
                (\ref{eq-3-6-connected-5:2})
                \\ \hline
                (\ref{eq-3-6-final-stage-3:4})
                & 
                (\ref{eq-3-6-third-stage-con-z-0}), 
                the RHS of (\ref{eq-3-6-3-stage-opera:4}), 
                (\ref{eq-3-6-3-stage-opera:6}), the RHS of (\ref{eq-3-6-3-stage-opera:5})
                \\ \hline
                (\ref{eq-3-6-final-stage-3-5})
                & 
                (\ref{eq-3-6-connected-3}) with (\ref{eq-3-6-connected-3:1}) substituted by $\bm{z}_{\rm o} \in \mathcal{Z}(\phi_2, \dddot{\bm{\alpha}}, \bm{\alpha}_+ \!-\! \bm{\alpha}_-  )$
               \\  \hline\hline
    \end{tabular*} 
     }
    \label{table-3-6-0}  
\end{table}



\section{Tractable Reformulation}

The TSDR-OTS model, with a complicated three-stage bi-level structure, is unable to be solved directly. Thus this section derives a traceable reformulation of the TSDR-OTS model.

\subsection{Reformulation of Recourse Functions}\label{sec-3-6-single-level}

We first reformulate the recourse functions into a more tractable form.

\subsubsection{Single-Level Reformulation}
The recourse function $Q_3$ is a bi-level optimization problem. Observing that the lower-level problem is a linear program, it can be replaced by its necessary and sufficient Karush–Kuhn–Tucker (KKT) conditions \cite{4-1279}:
\begin{subequations}\label{eq-3-6-KKT}
    \begin{align}
        \bm{W} \bm{x}_{\bm{o}}^{\bm{\varepsilon}} \leq \bm{u} - \bm{V} \![\bm{x}_{\bm{o}}^{\bm{\varepsilon}}]_{\rm z}, \bm{\lambda}_{\bm{o}}^{\bm{\varepsilon}} \geq \bm{0}, \bm{W}^T \bm{\lambda}_{\bm{o}}^{\bm{\varepsilon}} = \bm{l}  \label{eq-3-6-KKT:1} \\
        (\bm{\lambda}_{\bm{o}}^{\bm{\varepsilon}})^T ( \bm{u} - \bm{V} [\bm{x}_{\bm{o}}^{\bm{\varepsilon}}]_{\rm z} - \bm{W} [\bm{x}_{\bm{o}}^{\bm{\varepsilon}}]_{\rm c}  ) = 0 \label{eq-3-6-KKT:2}
    \end{align}
\end{subequations}
to obtain a single-level reformulation of $Q_3$. Here $\bm{\lambda}_{\bm{o}}^{\bm{\varepsilon}} \!\in\! \mathbb{R}^{n_{\rm u}}$ with $n_{\rm u}$ being the dimension of $\bm{u}$. The KKT complementarity conditions (\ref{eq-3-6-KKT:2}) can be further replaced by the mixed-integer reformulation 
\begin{equation}\label{eq-3-6-KKT-MI}
    \bm{u} - \bm{V} [\bm{x}_{\bm{\varepsilon}}]_{\rm z} - \bm{W} [\bm{x}_{\bm{\varepsilon}}]_{\rm c}  \leq M_{\rm p} (\bm{1} - \bm{\xi}_{\bm{o}}^{\bm{\varepsilon}} ), \bm{\lambda}_{\bm{o}}^{\bm{\varepsilon}} \leq M_{\rm d} \bm{\xi}_{\bm{o}}^{\bm{\varepsilon}}
\end{equation}
with $\bm{\xi}_{\bm{o}}^{\bm{\varepsilon}} \in \mathbb{B}^{n_{\rm u}}$ and sufficiently large constants $M_{\rm p}$ and $M_{\rm d}$.

\subsubsection{Relatively Complete Recourse (RCR) Reformulation}

The TSDR-OTS model has RCR iff the subproblem in each stage is feasible for any feasible decision made in the previous stages and for every realization of uncertainties in their supports. 
Generally, RCR of the TSDR-OTS model cannot always be ensured. 
For example, large VRE forecast errors following an inappropriate first-stage scheduling solution may cause incapability of the second-stage corrective control, i.e., infeasibility of the second-stage subproblem. 
Therefore, solving the TSDR-OTS model can be difficult since it is NP-hard just to find a feasible first-stage solution when RCR is not ensured \cite{4-1277}. To avoid this, we reformulate (\ref{eq-3-6-final-stage-2}) and (\ref{eq-3-6-final-stage-3}) as penalized problems. Specifically, the second-stage recourse function is reformulated as
\begin{subequations}\label{eq-3-6-reform-2}
    \begin{align}
        \!\!\!\! Q_2( \bm{x}, \bm{\varepsilon} ) \!:=\! 
        \min_{\tilde{\bm{x}}_{\bm{\varepsilon}} \in \mathbb{R}^{\tilde{n}_2}} ~ \!\!& \tilde{\bm{c}}_{\rm 2}^T \tilde{\bm{x}}_{\bm{\varepsilon}} 
        \!+\! \sup_{ O \in \mathcal{O} }   \mathbb{E}_{ O }   
        [ Q_3( \bm{x}, \bm{x}_{\bm{\varepsilon}}, \bm{o} ) ]   \label{eq-3-6-reform-2:1} \\
        \text{s.t.} ~ & 
        \tilde{\bm{C}} \tilde{\bm{x}}_{\bm{\varepsilon}} + \tilde{\bm{D}} \bm{x} \geq  \tilde{\bm{d}}( \bm{\varepsilon} )   \label{eq-3-6-reform-2:2} 
    \end{align}
\end{subequations}
where $\tilde{\bm{x}}_{\bm{\varepsilon}}$ is the augmentation of $\bm{x}_{\bm{\varepsilon}}$ with $y_{\bm{\varepsilon}} \!\in\! \mathbb{R}$, (\ref{eq-3-6-reform-2:2}) denotes (\ref{eq-3-6-final-stage-2}) with (\ref{eq-3-6-second-stage-cst:3}) slacked as
\begin{equation}
    \ddot{\bm{p}}_{\rm v}^{\rm max} \!\geq\! \dot{\bm{p}}_{\rm v}^{\rm max} \!+\! \bm{\varepsilon} \!-\! y_{\bm{\varepsilon}} \bm{1}, 
    \ddot{\bm{p}}_{\rm v}^{\rm max} \!\leq\! \dot{\bm{p}}_{\rm v}^{\rm max} \!+\! \bm{\varepsilon} \!+\! y_{\bm{\varepsilon}} \bm{1}, y_{\bm{\varepsilon}} \!\geq\! 0
\end{equation}
and $\tilde{\bm{c}}_{\rm 2}^T \tilde{\bm{x}}_{\bm{\varepsilon}} \!=\! \bm{c}_2^T \bm{x}_{\bm{\varepsilon}} + \sigma_2 y_{\bm{\varepsilon}}$ with $\sigma_2$ being the penalty coefficient.

Similarly, the penalized reformulation of (\ref{eq-3-6-final-stage-3}) can be formed by slacking constraints (\ref{eq-3-6-final-stage-3:2})-(\ref{eq-3-6-final-stage-3:4}), i.e., adding $\!-y_{\bm{o}}^{\bm{\varepsilon}} \bm{1}$ to the RHS of them, where $y_{\bm{o}}^{\bm{\varepsilon}} \!\in\! \mathbb{R}$ is the slack variable; and replacing the objective function in (\ref{eq-3-6-final-stage-3:1}) by $\tilde{\bm{c}}_{\rm 3}^T \tilde{\bm{x}}_{\bm{o}}^{\bm{\varepsilon}} \!=\! \bm{c}_3^T \bm{x}_{\bm{o}}^{\bm{\varepsilon}} \!+ \sigma_3 y_{\bm{o}}^{\bm{\varepsilon}} $, with $\sigma_3$ being the penalty coefficient and $\tilde{\bm{x}}_{\bm{o}}^{\bm{\varepsilon}} = [(\bm{x}_{\bm{o}}^{\bm{\varepsilon}})^T ~ \bm{\lambda}^T ~ \bm{\xi}^T ~ y_{\bm{o}}^{\bm{\varepsilon}}]^T$. Note that (\ref{eq-3-6-final-stage-3-5}) is always feasible and thus remains unchanged in the penalized reformulation. 
Then, according to the single-level reformulation, substituting (\ref{eq-3-6-final-stage-3-5}) in the penalized reformulation by (\ref{eq-3-6-KKT:1}) and (\ref{eq-3-6-KKT-MI}) gives
\begin{subequations}\label{eq-3-6-reform}
    \begin{align}
        \!\!\!\! Q_3(\! \bm{x},\! \bm{x}_{\bm{\varepsilon}},\! \bm{o} \!)   &  \!:=\! 
        \min\nolimits_{\tilde{\bm{x}}_{\bm{o}}^{\bm{\varepsilon}} \in \mathbb{R}^{\tilde{n}_3} \!\times \mathbb{B}^{\tilde{m}_3}} ~   \tilde{\bm{c}}_{\rm 3}^T \tilde{\bm{x}}_{\bm{o}}^{\bm{\varepsilon}}  \label{eq-3-6-reform:1}\\[-1mm] 
        \text{s.t.} ~&  \tilde{\bm{F}} \tilde{\bm{x}}_{\bm{o}}^{\bm{\varepsilon}} \!+\! \tilde{\bm{H}} {\bm{x}}_{\bm{\varepsilon}} \!+\!  \tilde{\bm{L}} \bm{x} \!\geq\! \tilde{\bm{f}} \label{eq-3-6-reform:2} \\[-1mm]  
        \!\!\!~ &  \tilde{\bm{M}} \tilde{\bm{x}}_{\bm{o}}^{\bm{\varepsilon}} + \tilde{\bm{N}} \bm{o} \geq \tilde{\bm{g}}  \label{eq-3-6-reform:3} \\[-1mm]  
        \!\!\!~ & \tilde{\bm{R}} \tilde{\bm{x}}_{\bm{o}}^{\bm{\varepsilon}} \!+\! \tilde{\bm{S}}(\bm{o}) [\bm{x}^T ~\bm{x}_{\bm{\varepsilon}}^T ~(\bm{x}_{\bm{o}}^{\bm{\varepsilon}})^T ]^T  \!\geq\! \tilde{\bm{h}}   \label{eq-3-6-reform:4}
    \end{align}
\end{subequations}
where 
(\ref{eq-3-6-reform:2}) represents constraints (\ref{eq-3-6-final-stage-3:2}) after being slacked, (\ref{eq-3-6-KKT:1}), and (\ref{eq-3-6-KKT-MI}); 
(\ref{eq-3-6-reform:3}) and (\ref{eq-3-6-reform:4}) represent (\ref{eq-3-6-final-stage-3:3}) and (\ref{eq-3-6-final-stage-3:4}) both after being slacked, respectively. Thereby, the TSDR-OTS model always has RCR provided that (\ref{eq-3-6-final-stage-1:2}) is feasible.
It is also noted that practical systems generally have adequate control resources responding to VRE uncertainty and contingencies, and thus the TSDR-OTS model (\ref{eq-3-6-final-stage-1}) is assumed to be feasible. In this case, the above penalized reformulations will not cause infeasibility of the final solution to (\ref{eq-3-6-final-stage-1}) since all slack variables can be penalized to zero at the optimum.


\subsection{Reformulation of the TSDR-OTS Model}

Next, we further derive the tractable TSDR-OTS model. Denote $\sup\nolimits_{O \in \mathcal{O}}  \mathbb{E}_{O} [Q_3( \bm{x}, \bm{x}_{\bm{\varepsilon}}, \bm{o} )] $ by $\mathcal{Q}_3( \bm{x}, \bm{x}_{\bm{\varepsilon}})$. We have
\begin{subequations}\label{eq-3-6-before-duality}
    \begin{align}
        & \mathcal{Q}_3(\bm{x}, \bm{x}_{\bm{\varepsilon}}) = \sup_{ O \in \mathcal{O} }   \sum\nolimits_{\bm{o} \in \Omega_{\bm{\varepsilon}} }  Q_3( \bm{x}, \bm{x}_{\bm{\varepsilon}}, \bm{o} ) O( \bm{o}) \label{eq-3-6-before-duality:1}\\[-1mm]
        & \text{s.t.}  \sum\nolimits_{\bm{o} \in \Omega_{\bm{\varepsilon}}}\!\! \bm{T} (\bm{1} - \bm{o}) O(\bm{o}) \leq \tilde{\bm{o}}, \sum\nolimits_{\bm{o} \in \Omega_{\bm{\varepsilon}}} \!\!  O(\bm{o}) = 1  \label{eq-3-6-before-duality:2}
    \end{align}
\end{subequations}
Here $\Omega_{\bm{\varepsilon}} = \Omega$ while $\Omega_{\bm{\varepsilon}}$ will be specified differently later. By strong duality \cite{4-1279}, (\ref{eq-3-6-before-duality}) admits the following dual formulation:
\begin{subequations}\label{eq-3-6-duality}
    \begin{align}
        \!\!\!\!  \mathcal{Q}_3(\bm{x}, \bm{x}_{\bm{\varepsilon}}) \!\!=\!\! ~& \min_{ \bm{\lambda}_{\bm{\varepsilon}}, \lambda_{\bm{\varepsilon}}' }  \tilde{\bm{o}}^T \bm{\lambda}_{\bm{\varepsilon}}  + \lambda_{\bm{\varepsilon}}' \label{eq-3-6-duality:1}\\
        \text{s.t.} ~&  \bm{\lambda}_{\bm{\varepsilon}}^T \bm{T} (\bm{1} \!-\! \bm{o}) \!+\! \lambda_{\bm{\varepsilon}}' \!\geq\! Q_3( \bm{x}, \bm{x}_{\bm{\varepsilon}}, \bm{o} ) ~ \forall \bm{o} \!\in\! \Omega_{\bm{\varepsilon}} \label{eq-3-6-duality:2}\\
        & \bm{\lambda}_{\bm{\varepsilon}} \geq \bm{0} \label{eq-3-6-duality:3}
    \end{align}
\end{subequations}
where $\bm{\lambda}_{\bm{\varepsilon}}$ and $\lambda_{\bm{\varepsilon}}'$ are dual variables associated with the two constraints of (\ref{eq-3-6-before-duality:2}), respectively. We further replace $Q_3( \bm{x}, \bm{x}_{\bm{\varepsilon}}, \bm{o} )$ in (\ref{eq-3-6-duality:1}) by the objective function in (\ref{eq-3-6-reform}) and augment constraints of (\ref{eq-3-6-duality}) with that of (\ref{eq-3-6-reform}), which yields
\begin{subequations}\label{eq-3-6-duality-sub}
    \begin{align}
        \!\!\!\!\!\! \mathcal{Q}_3( \bm{x}, \bm{x}_{\bm{\varepsilon}}) & \!\!=\!\!\!  \min_{ \bm{\lambda}_{\bm{\varepsilon}}, \lambda_{\bm{\varepsilon}}', \{ \tilde{\bm{x}}_{\bm{o}}^{\bm{\varepsilon}} \} } ~ \tilde{\bm{o}}^T \bm{\lambda}_{\bm{\varepsilon}}  + \lambda_{\bm{\varepsilon}}' \label{eq-3-6-duality-sub:1}\\
        \text{s.t.} ~&  \bm{\lambda}_{\bm{\varepsilon}}^T \bm{T} (\bm{1} \!\!-\!\! \bm{o}) \!\!+\!\! \lambda_{\bm{\varepsilon}}' \!\!\geq\!\!  \tilde{\bm{c}}_{\rm 3}^T \tilde{\bm{x}}_{\bm{o}}^{\bm{\varepsilon}}, \text{(\ref{eq-3-6-reform:3}), (\ref{eq-3-6-reform:4})} ~ \forall \bm{o} \!\in\! \Omega_{\bm{\varepsilon}} \label{eq-3-6-duality-sub:2}\\
        & \text{(\ref{eq-3-6-duality:3}), (\ref{eq-3-6-reform:2}) } \!\!\!\!\!\! \label{eq-3-6-duality-sub:3} 
    \end{align}
\end{subequations}
The equality in (\ref{eq-3-6-duality-sub:1}) holds since projection of the feasible region of the optimization problem in (\ref{eq-3-6-duality-sub}) on the $\bm{\lambda}_{\bm{\varepsilon}}$-$\lambda_{\bm{\varepsilon}}'$ plane equals to the feasible region of the optimization problem in (\ref{eq-3-6-duality}), and objective functions of these two optimization problems both only depend on $\bm{\lambda}_{\bm{\varepsilon}}$ and $\lambda_{\bm{\varepsilon}}'$.
 
Substituting (\ref{eq-3-6-duality-sub}) into (\ref{eq-3-6-reform-2}) gives the equivalent second-stage recourse function
\begin{equation}\label{eq-3-6-stage2-final} 
    \begin{aligned}
        Q_2( \bm{x}, \bm{\varepsilon} ) :=
        \min_{ \tilde{\bm{x}}_{\bm{\varepsilon}}, \{\tilde{\bm{x}}_{\bm{o}}^{\bm{\varepsilon}} \}, \bm{\lambda}_{\bm{\varepsilon}}, \lambda_{\bm{\varepsilon}}' }  & \tilde{\bm{c}}_2^T \tilde{\bm{x}}_{\bm{\varepsilon}} + \tilde{\bm{o}}^T \!\bm{\lambda}_{\bm{\varepsilon}}  + \lambda_{\bm{\varepsilon}}' \\
          \text{s.t.} & 
        \text{ (\ref{eq-3-6-reform-2:2}), (\ref{eq-3-6-duality-sub:2}), (\ref{eq-3-6-duality-sub:3}) } 
    \end{aligned}
\end{equation}
We further replace the expectation in (\ref{eq-3-6-final-stage-1}) by its scenario-based approximation, yielding the final tractable reformulation of the TSDR-OTS model as follows:
\begin{equation}\label{eq-3-6-final-opt-model}
    \begin{aligned}
        & Q_1 \!=\!\! \min_{\bm{x}, \{\! \tilde{\bm{x}}_{\bm{\varepsilon}} \!\}, \{\! \tilde{\bm{x}}_{\bm{o}}^{\bm{\varepsilon}} \!\}, \{\!\bm{\lambda}_{\bm{\varepsilon}}\!\}, \{\!\lambda_{\bm{\varepsilon}}' \!\} } 
        \bm{c}_{\rm 1}^T \bm{x} \!\!+\!\!\! \sum_{ \bm{\varepsilon} \in \Xi }  \!\! {P}(\! \bm{\varepsilon} \!)\! \left[ \tilde{\bm{c}}_2^T \tilde{\bm{x}}_{\bm{\varepsilon}} \!\!+\!\! \tilde{\bm{o}}^T \bm{\lambda}_{\bm{\varepsilon}}  \!\!+\!\! \lambda_{\bm{\varepsilon}}' \right]  \\
        &~~~~~~~~~~~~~~~~~~~ \text{s.t.}   \text{(\ref{eq-3-6-final-stage-1:2})}, \{ \text{(\ref{eq-3-6-reform-2:2}), (\ref{eq-3-6-duality-sub:2}), (\ref{eq-3-6-duality-sub:3}) } | \forall \bm{\varepsilon} \!\in\! \Xi  \}
    \end{aligned}
\end{equation}
where $\Xi$ is the set of finite number of scenarios of $\bm{\varepsilon}$.

\section{Solution Approach}

To solve the tractable reformulation (\ref{eq-3-6-final-opt-model}) efficiently, in this section, we develop a solution approach which combines the DW procedure with the nested CCG algorithm. In particular, for the inner CCG loop, we utilize strong duality to derive tractable formulation of the multi-level master problem (MP), preventing introducing massive auxiliary binary variables.



\vspace{-4pt}
\subsection{Overall Framework of the Solution Approach}

\begin{algorithm}
    \nonl \includegraphics[width=0.955\linewidth]{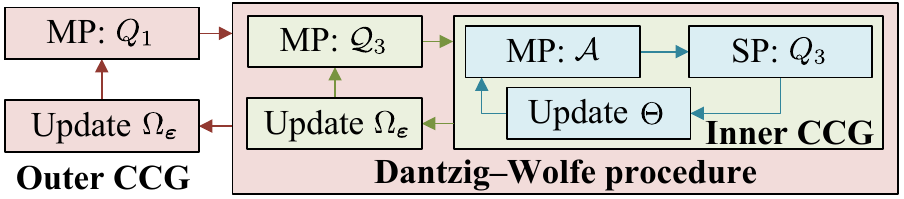} 

    \vspace{-9pt}

    \nonl \noindent\rule{8.5cm}{0.5pt}

    \DontPrintSemicolon
\small{ 
    \KwInput{$\epsilon_1$, $\epsilon_2$, $\epsilon_3$, $n_{\max}$, $n_{\rm o}$}
    \KwOutput{$\bm{x}^{(i)}$} 

    Initialization: $i\!=\!0$, $Q_{\rm lb}^{(i)} \!=\! -\infty$, $Q_{\rm ub}^{(i)} \!=\! \infty$; $\forall \bm{\varepsilon} \in \Xi$, $\Omega_{\bm{\varepsilon}}^{(i)} \!=\! \emptyset$

    \While{$(Q_{\rm ub}^{(i)} - Q_{\rm lb}^{(i)})/Q_{\rm ub}^{(i)} \!\leq\! \epsilon_1 $}
    {
        $(\bm{x}^{(i)}\!,\! \bm{x}_{\bm{\varepsilon}}^{(i)}\!,\! \cdots \!) \overset{\text{arg}}{\gets} Q_1(\! \{\Omega_{\bm{\varepsilon}}^{(i)} \!| \forall \bm{\varepsilon} \!\in\! \Xi \} \!)$ 

        $Q_{\rm lb}^{(i+1)} \gets Q_1( \{\Omega_{\bm{\varepsilon}}^{(i)} | \forall \bm{\varepsilon} \in \Xi \} )$

        \If{$(Q_{\rm ub}^{(i)} - Q_{\rm lb}^{(i+1)})/Q_{\rm ub}^{(i)} \!\leq\! \epsilon_1 $}
        {
            \textbf{Break} 
        }

        \For{$\bm{\varepsilon} \in \Xi $}    
        { 
            Initialization: $j_{\bm{\varepsilon}} = 0$, $\Omega_{\bm{\varepsilon}}^{(i, j_{\bm{\varepsilon}})} \gets \Omega_{\bm{\varepsilon}}^{(i)}$

            \Repeat{$\mathcal{A}(\bm{x}^{(i)}, \bm{x}_{\bm{\varepsilon}}^{(i)}\!,\! \bm{\lambda}_{\bm{\varepsilon}}^{(j_{\bm{\varepsilon}}-1)}\!,\! \lambda_{\bm{\varepsilon}}'^{(j_{\bm{\varepsilon}}-1)} ) \!\!\leq\!\! \epsilon_2 $ {\rm or} $j_{\bm{\varepsilon}} \!\!>\!\! n_{\max}$ }
            {

                $(\bm{\lambda}_{\bm{\varepsilon}}^{(\!j_{\!\bm{\varepsilon}}\!)}\!,\! \lambda_{\bm{\varepsilon}}'^{(\!j_{\!\bm{\varepsilon}}\!)}\!,\! \{\! \tilde{\bm{x}}_{\bm{o}}^{\bm{\varepsilon} (\!j_{\bm{\varepsilon}}\!) } \!\}\!,\! O^{(\!j_{\bm{\varepsilon}}\!)}\!,\!...\!) \!\overset{\text{arg}}{\gets}\! \mathcal{Q}_3(\bm{x}^{(i)}\!,\! \bm{x}_{\bm{\varepsilon}}^{(i)}\!,\!  \Omega_{\bm{\varepsilon}}^{(\!i,j_{\bm{\varepsilon}}\!)} \!)$ 


                Initialization: $k\!=\!0$, $-\mathcal{A}_{\rm lb}^{(k)} \!=\! \mathcal{A}_{\rm ub}^{(k)} \!=\! \infty$, $\Theta^{(k)} \!\!=\! \emptyset$

                \While{$(\mathcal{A}_{\rm ub}^{(k)} \!-\! \mathcal{A}_{\rm lb}^{(k)})/\mathcal{A}_{\rm ub}^{(k)} \!\leq\! \epsilon_3 $}
                {
                    $(\bm{o}^{(j_{\bm{\varepsilon}}, k)}\!,\! \cdots )  \overset{\text{arg}}{\gets} \mathcal{A}(\bm{x}^{(i)}\!,\! \bm{x}_{\bm{\varepsilon}}^{(i)}\!,\! \bm{\lambda}_{\bm{\varepsilon}}^{(j_{\bm{\varepsilon}})}\!,\!  \lambda_{\bm{\varepsilon}}'^{(j_{\bm{\varepsilon}})}\!,\! \Theta^{(k)})$ 

                    $\mathcal{A}_{\rm ub}^{(k+1)} \gets \mathcal{A}(\bm{x}^{(i)},\! \bm{x}_{\bm{\varepsilon}}^{(i)},\! \bm{\lambda}_{\bm{\varepsilon}}^{(j_{\bm{\varepsilon}})},\!  \lambda_{\bm{\varepsilon}}'^{(j_{\bm{\varepsilon}})}, \Theta^{(k)}) $

                    \If{$(\mathcal{A}_{\rm ub}^{(k+1)} \!-\! \mathcal{A}_{\rm lb}^{(k)})/\mathcal{A}_{\rm ub}^{(k+1)} \!\leq\! \epsilon_2 $}{\textbf{Break}}

                    $(\bm{y}_{\bm{o}}^{\bm{\varepsilon} (k) }, \cdots ) \overset{\text{arg}}{\gets} Q_3(\bm{x}^{(i)}, \bm{x}_{\bm{\varepsilon}}^{(i)}, \bm{o}^{(j_{\bm{\varepsilon}}, k)})$ 

                    $\Theta^{(k+1)} \gets \Theta^{(k)} \cup \{ \bm{y}_{\bm{o}}^{\bm{\varepsilon} (k) } \} $

                    $\mathcal{A}_{\rm lb}^{(k\!+\!1)} \!\!\gets\! 
                    \max\!\! \left\{\!\!
                        \begin{aligned}
                            & \mathcal{A}_{\rm lb}^{(k)}\!,\! Q_3( \bm{x}^{(i)}\!, \bm{x}_{\bm{\varepsilon}}^{(i)}\!,\! \bm{o}^{(j_{\bm{\varepsilon}}, k)} \!) -\! \\
                            &[(\!\bm{\lambda}_{\bm{\varepsilon}}^{(j_{\bm{\varepsilon}})}\!)\!^T \!\bm{T} \!(\!\bm{1} \!\!-\!\! \bm{o}^{(j_{\bm{\varepsilon}}, k)} \!)\! \!+\!\! \lambda_{\bm{\varepsilon}}'^{(j_{\bm{\varepsilon}})} ] 
                        \end{aligned}
                        \! \right\}  
                    $
                
                    $k = k + 1$
                }

                $\Omega_{\bm{\varepsilon}}^{(i, j_{\bm{\varepsilon}} + 1)} \gets \Omega_{\bm{\varepsilon}}^{(i, j_{\bm{\varepsilon}})} \cup \{\bm{o}^{(j_{\bm{\varepsilon}}, k-1)} \} $, 
                $j_{\bm{\varepsilon}} \gets j_{\bm{\varepsilon}} + 1$

            }

            $\Omega_{\bm{\varepsilon}}^{(\!i\!+\!1\!)} \!\!\gets\! \mathcal{C}\bm{(} \Omega_{\bm{\varepsilon}}^{(\!i, j_{\bm{\varepsilon}}\!-\!1\!)} \!\backslash\! \Omega_{\bm{\varepsilon}}^{(i)}\!,\! \tilde{\bm{c}}_{\rm 3}^T \tilde{\bm{x}}_{\bm{o}}^{\bm{\varepsilon} (j_{\bm{\varepsilon}}\!-\!1) } \!O^{(\!j_{\bm{\varepsilon}}\!-\!1\!)}\!(\!\bm{o}\!), n_{\rm o} \bm{)} \!\cup\! \Omega_{\bm{\varepsilon}}^{(\!i\!)} $

        }

        $Q_{\rm ub}^{(i+1)} \!\!\gets\! \bm{c}_{\rm 1}^T \bm{x}^{(i)} \!\!+\!\!\! \sum\limits_{ \bm{\varepsilon} \in \Xi }  \!\! {P}(\! \bm{\varepsilon} \!)\!\! 
        \left[\!
            \begin{aligned}
                & \tilde{\bm{c}}_2^T \!\tilde{\bm{x}}_{\bm{\varepsilon}}^{(i)} \!\!+\!\!   \mathcal{Q}_3(\!\bm{x}^{\!(i)}\!, \bm{x}_{\bm{\varepsilon}}^{(i)}\!, \Omega_{\bm{\varepsilon}}^{(i,j_{\bm{\varepsilon}}\!-\!1)} \!) +\! \\[-1mm]
                & \mathcal{A}(\bm{x}^{(i)}, \bm{x}_{\bm{\varepsilon}}^{(i)}, \bm{\lambda}_{\bm{\varepsilon}}^{(j_{\bm{\varepsilon}}-1)}, \lambda_{\bm{\varepsilon}}'^{(j_{\bm{\varepsilon}}-1)} )   
            \end{aligned}  
        \right] \!\!\!$
        
        $i \gets i + 1$

    }
    \caption{The proposed solution approach}
    \label{alg-3-6-1}
}

\justifying
\footnotesize{
\nonl \noindent  Notation: $Q_1( \{\Omega_{\bm{\varepsilon}}' | \forall \bm{\varepsilon} \!\in\! \Xi \} )$, $Q_2( \bm{x}, \bm{\varepsilon}, \Omega_{\bm{\varepsilon}}' )$, $\mathcal{Q}_3(\bm{x}, \bm{x}_{\bm{\varepsilon}},  \Omega_{\bm{\varepsilon}}' )$ and $\mathcal{A}(\bm{x},\! \bm{x}_{\bm{\varepsilon}},\! \bm{\lambda}_{\bm{\varepsilon}},\!  \lambda_{\bm{\varepsilon}}', \Theta') $ denote $Q_1$, $Q_2( \bm{x}, \bm{\varepsilon})$ and $\mathcal{Q}_3(\bm{x}, \bm{x}_{\bm{\varepsilon}} )$ all with $\Omega_{\bm{\varepsilon}} \!\!=\!\! \Omega_{\bm{\varepsilon}}'$, and $\mathcal{A}(\bm{x},\! \bm{x}_{\bm{\varepsilon}},\! \bm{\lambda}_{\bm{\varepsilon}},\! \lambda_{\bm{\varepsilon}}' )$ with $\Theta \!\!=\!\! \Theta'$, respectively; $x \!\overset{\text{arg}}{\gets}\!\! X$ represents solving problem $X$ and assigning the optimal solution to $x$, $X$ the objective function value at the optimum, and $\mathcal{C}(Y, y , n)$ the $n$ elements in set $Y$ whose associated values of $y$ are the largest. The superscript ``$(i)$'' indicates association with the $i$-th iteration. 
}
\end{algorithm}

Algorithm \ref{alg-3-6-1} gives the main framework and pseudocode of the proposed solution approach. In the outer CCG loop, for each iteration $i$, the MP $Q_1( \{\Omega_{\bm{\varepsilon}}^{(i)} | \forall \bm{\varepsilon} \!\in\! \Xi \} )$ with $\Omega_{\bm{\varepsilon}}^{(i)} \!\subset\! \Omega$ being a subset of all $N\!-\!k$ contingencies, is firstly solved to obtain a trial solution $(\bm{x}^{(i)}, \bm{x}_{\bm{\varepsilon}}^{(i)}, \cdots)$ and its associated lower bound of $Q_1$, i.e., $Q_{\rm lb}^{(i+1)}$. Then the DW procedure is executed for each $\bm{\varepsilon} \!\in\! \Xi$ to obtain an upper bound of $Q_1$, i.e., $Q_{\rm ub}^{(i+1)}$, and update $\Omega_{\bm{\varepsilon}}^{(i)}$ to $\Omega_{\bm{\varepsilon}}^{(i+1)}$ which determines the new constraints added to the MP in the next iteration. The inner CCG loop is designed to solve the subproblems (SPs) in the DW procedure. 
The algorithm terminates when the relative gap between the upper and lower bounds, i.e., $(Q_{\rm ub}^{(i)} \!-\! Q_{\rm lb}^{(i)})/Q_{\rm ub}^{(i)}$, is below tolerance $\epsilon_1$. Finally, the $\epsilon_1$-optimal solution $\bm{x}^{(i)}$, given by the last iteration, is outputted for the dispatch.

It is noted that computation time of the DW procedure from line 7 to 23 in Algorithm \ref{alg-3-6-1} takes up the majority of the whole, which however, can be significantly reduced by executing the DW procedure of each $\bm{\varepsilon} \in \Xi$ or each group of $\bm{\varepsilon} \in \Xi$ in different computing units in parallel.

\vspace{-4pt}
\subsection{Dantzig–Wolfe Procedure}

In iteration $i$ of the outer CCG loop, for each $\bm{\varepsilon} \in \Xi$, the counter $j_{\bm{\varepsilon}}$ and contingency subset $\Omega_{\bm{\varepsilon}}^{(\!i,j_{\bm{\varepsilon}}\!)}$ are initialized to 0 and $\Omega_{\bm{\varepsilon}}^{(i)}$ respectively, and the DW procedure performs up to $j_{\bm{\varepsilon}} = n_{\max}$ iterations. In each iteration $j_{\bm{\varepsilon}}$ of the DW procedure, the third-stage recourse function problem which considers only contingencies in $\Omega_{\bm{\varepsilon}}^{(\!i,j_{\bm{\varepsilon}}\!)}$, i.e.,  $\mathcal{Q}_3(\bm{x}^{(i)}, \bm{x}_{\bm{\varepsilon}}^{(i)},  \Omega_{\bm{\varepsilon}}^{(i,j_{\bm{\varepsilon}})} )$, is firstly solved to update the variables associate with this recourse function, i.e., 
$(\bm{\lambda}_{\bm{\varepsilon}}^{(j_{\bm{\varepsilon}})}, \lambda_{\bm{\varepsilon}}'^{(j_{\bm{\varepsilon}})}, \{ \tilde{\bm{x}}_{\bm{o}}^{\bm{\varepsilon} (j_{\bm{\varepsilon}}) } \}, O^{(j_{\bm{\varepsilon}})},...)$. In addition, $\mathcal{Q}_3(\bm{x}^{(i)}, \bm{x}_{\bm{\varepsilon}}^{(i)},  \Omega_{\bm{\varepsilon}}^{(i,j_{\bm{\varepsilon}})} )$ gives a lower bound of $\mathcal{Q}_3(\bm{x}^{(i)}, \bm{x}_{\bm{\varepsilon}}^{(i)} )$. 
Next, $\bm{o} \!\in\! \Omega$ with the highest reduction of the third-stage objective function is identified by solving 
$\mathcal{A}(\bm{x}^{(i)},\! \bm{x}_{\bm{\varepsilon}}^{(i)},\! \bm{\lambda}_{\bm{\varepsilon}}^{(j_{\bm{\varepsilon}})},\!  \lambda_{\bm{\varepsilon}}'^{(j_{\bm{\varepsilon}})} )$ which is defined as 
\vspace{-4pt}
\begin{equation}\label{eq-3-6-problem-A}
    \vspace{-4pt}
    \mathcal{A}(\bm{x}, \bm{x}_{\bm{\varepsilon}}, \bm{\lambda}_{\bm{\varepsilon}},  \lambda_{\bm{\varepsilon}}' ) \!:=\! \max_{\bm{o} \in \Omega }~  Q_3( \bm{x},\! \bm{x}_{\bm{\varepsilon}},\! \bm{o} ) \!-\!  [\bm{\lambda}_{\bm{\varepsilon}}^T \bm{T} (\bm{1} \!-\! \bm{o}) \!+\! \lambda_{\bm{\varepsilon}}' ]
\end{equation}
The identified contingency $\bm{o}$ is added to $\Omega_{\bm{\varepsilon}}^{(i,j_{\bm{\varepsilon}})}$ for the next iteration. The DW procedure terminates when the highest reduction of the third-stage objective function is below the given tolerance $\epsilon_2$ or the iteration number reaches $n_{\max}$. Among all identified new contingencies, the first $n_{\rm o}$ ones with the largest contribution to the objective value of $\mathcal{Q}_3(\bm{x}^{(i)}, \bm{x}_{\bm{\varepsilon}}^{(i)},  \Omega_{\bm{\varepsilon}}^{(i,j_{\bm{\varepsilon}})} )$ in the last iteration of the DW procedure are added to $\Omega_{\bm{\varepsilon}}^{(i)}$. When the DW procedure is performed for all $\bm{\varepsilon} \!\in\! \Xi$, a new upper bound of $Q_1$ is computed as shown in line 24 of Algorithm \ref{alg-3-6-1}.$\!$


 
\vspace{-4pt}
\subsection{Inner CCG Loop}

The solution of $\mathcal{A}$ in the DW procedure is obtained through the inner CCG loop, as shown in lines 12 to 20 in Algorithm \ref{alg-3-6-1}. The formulation of $\mathcal{A}$ utilized by the inner CCG loop is derived below. 
We first reformulate (\ref{eq-3-6-reform:4}) by the Big-M method as
\begin{equation}
    \begin{bmatrix}
        \!\tilde{\bm{R}} \tilde{\bm{x}}_{\bm{o}}^{\bm{\varepsilon}} \!\!+\!\! \tilde{\bm{S}}(\!\bm{1}_N\!) [\bm{x}^T ~\!\bm{x}_{\bm{\varepsilon}}^T ~\!(\!\bm{x}_{\bm{o}}^{\bm{\varepsilon}}\!)^T ]^T 
\!\!\!+\!\! [ \tilde{\bm{S}}(\!\bm{1}_N\!) \bm{1} \!\!-\!\! \tilde{\bm{S}}(\!\bm{o}\!) \bm{1} ] \!M\!\! \\
\tilde{\bm{R}} \tilde{\bm{x}}_{\bm{o}}^{\bm{\varepsilon}} + [ \tilde{\bm{S}}(\bm{o}) \bm{1} ] M
    \end{bmatrix}
    \!\!\!\geq\!\!\! 
    \begin{bmatrix}
        \tilde{\bm{h}}  \\
        \tilde{\bm{h}}  
    \end{bmatrix}
\end{equation}
Divide $\tilde{\bm{x}}_{\bm{o}}^{\bm{\varepsilon}}$ into continuous and binary sub-vectors, denoted as $\bm{\gamma}_{\bm{o}}^{\bm{\varepsilon}} \in \mathbb{R}^{\tilde{n}_3}$ and ${\bm{y}}_{\bm{o}}^{\bm{\varepsilon}} \in \mathbb{B}^{\tilde{m}_3}$, respectively. 
Then (\ref{eq-3-6-reform}) can be rewritten as the following compact form:
\begin{subequations}\label{eq-3-6-reform-ccg}
    \begin{align}
        \!\!\!\! Q_3(\! \bm{x},\! \bm{x}_{\bm{\varepsilon}},\! \bm{o} \!)   &  \!:=\! 
        \min\nolimits_{ \bm{\gamma}_{\bm{o}}^{\bm{\varepsilon}} \in \mathbb{R}^{\tilde{n}_3}, {\bm{y}}_{\bm{o}}^{\bm{\varepsilon}} \in \mathbb{B}^{\tilde{m}_3}  } ~   \bm{a}^T  \bm{\gamma}_{\bm{o}}^{\bm{\varepsilon}} \label{eq-3-6-reform-ccg:1}\\
        \text{s.t.} ~&
        \bm{P} \bm{x} + \bm{Q} \bm{x}_{\bm{\varepsilon}} + \bm{U}_{\rm c}  \bm{\gamma}_{\bm{o}}^{\bm{\varepsilon}} + \bm{U}_{\rm b} {\bm{y}}_{\bm{o}}^{\bm{\varepsilon}} +  \bm{O} \bm{o} \geq \bm{t}
        \label{eq-3-6-reform-ccg:2}  
    \end{align}
\end{subequations}
Note that $Q_3$ originally does not satisfy the extended RCR property, i.e., (\ref{eq-3-6-reform}) is feasible with $\bm{x},\! \bm{x}_{\bm{\varepsilon}},\! \bm{o}$ and ${\bm{y}}_{\bm{o}}^{\bm{\varepsilon}}$ fixed to their any possible values. Similarly to the derivation of (\ref{eq-3-6-reform}) to satisfy the RCR property, we can add the term $y_{\bm{o}}^{\bm{\varepsilon}} \bm{1}$ to constraints (\ref{eq-3-6-KKT-MI}) to ensure the extended RCR property. We assume that this is already embodied in (\ref{eq-3-6-reform-ccg:2}). 

Substituting (\ref{eq-3-6-reform-ccg}) into (\ref{eq-3-6-problem-A}) yields the equivalent tri-level formulation as follows:
\begin{equation}\label{eq-3-6-function-A-tri-level}
    \begin{aligned}
        &\!\!\!\!\!\! \mathcal{A}(\bm{x}, \bm{x}_{\bm{\varepsilon}}, \bm{\lambda}_{\bm{\varepsilon}},  \lambda_{\bm{\varepsilon}}' ) \!:=\!  \\ 
        &\!\!\!\!\!\! \max_{\bm{o} \in \Omega } -\!  [\bm{\lambda}_{\bm{\varepsilon}}^T \bm{T} (\bm{1} \!-\! \bm{o}) \!\!+\!\! \lambda_{\bm{\varepsilon}}' ] \!\!+\!\!  \min_{ \bm{y}_{\bm{o}}^{\bm{\varepsilon}} \in \Theta }
        \! \{\! \min_{\bm{\gamma}_{\bm{o}}^{\bm{\varepsilon}} \in \mathbb{R}^{\tilde{n}_3}} \bm{a}^T  \bm{\gamma}_{\bm{o}}^{\bm{\varepsilon}}  ~\text{s.t. (\ref{eq-3-6-reform-ccg:2})} \! \}
    \end{aligned}
\end{equation}
with $\Theta = \mathbb{B}^{\tilde{m}_3}$. The innermost minimization problem is a linear program. Dualizing it gives
\begin{equation}\label{eq-3-6-lp-duality} 
        \!\!\!\! \max_{\bm{\pi}_{\bm{o}}^{\bm{\varepsilon}} \in \mathbb{R}^{n_{\rm t}} }    (\!\bm{t} \!-\!\bm{P} \bm{x} \!-\! \bm{Q} \bm{x}_{\bm{\varepsilon}} \!-\! \bm{U}_{\rm b} {\bm{y}}_{\bm{o}}^{\bm{\varepsilon}} \!-\!  \bm{O} \bm{o}\!)^T \!\bm{\pi}_{\bm{o}}^{\bm{\varepsilon}}  ~ \text{s.t.}~   \bm{\pi}_{\bm{o}}^{\bm{\varepsilon}} \!\!\geq\!\! \bm{0}, \bm{U}_{\rm c}^T \!\bm{\pi}_{\bm{o}}^{\bm{\varepsilon}} \!\!=\!\! \bm{a}
\end{equation}
with $n_{\rm t}$ being the dimension of vector $\bm{t}$, and $\bm{\pi}_{\bm{o}}^{\bm{\varepsilon}}$ the vector of dual variables.

By substituting (\ref{eq-3-6-lp-duality}) into (\ref{eq-3-6-function-A-tri-level}) and enumerating all possible assignments of $\bm{y}_{\bm{o}}^{\bm{\varepsilon}} \!\!\in\!\! \Theta$, we can obtain the monolithic equivalent form of  (\ref{eq-3-6-function-A-tri-level}) as follows:
\begin{subequations}\label{eq-3-6-function-A-equivalent-dual-emu}
    \begin{align}
        &\!\!\!\! \mathcal{A}(\bm{x},\! \bm{x}_{\bm{\varepsilon}},\! \bm{\lambda}_{\bm{\varepsilon}},\!  \lambda_{\bm{\varepsilon}}' ) \!\!:=\!\!\!\! \max_{\bm{o} \in \Omega, \varpi_{\bm{o}}^{\bm{\varepsilon}} \in \mathbb{R}, \bm{\pi}_{\!\bm{o}, \bm{y}}^{\bm{\varepsilon}} \!\in \mathbb{R}^{n_{\rm t}}  } \!\!\! \varpi_{\bm{o}}^{\bm{\varepsilon}} \!\!-\!\!  [\bm{\lambda}_{\bm{\varepsilon}}^T \bm{T} (\bm{1} \!\!-\!\! \bm{o}) \!\!+\!\! \lambda_{\bm{\varepsilon}}' ] \\
        &\!\!\!\! \text{s.t.}~   \varpi_{\bm{o}}^{\bm{\varepsilon}} \!\leq\! (\bm{t} \!-\!\bm{P} \bm{x} \!-\! \bm{Q} \bm{x}_{\bm{\varepsilon}} \!-\! \bm{U}_{\rm b} {\bm{y}}_{\bm{o}}^{\bm{\varepsilon}} \!-\!  \bm{O} \bm{o})^T \bm{\pi}_{\bm{o}, \bm{y} }^{\bm{\varepsilon}} ~ \forall \bm{y}_{\bm{o}}^{\bm{\varepsilon}} \!\in\! \Theta \label{eq-3-6-function-A-equivalent-dual-emu:1}\\
        &\!\!\!\! ~~~~  \bm{\pi}_{\bm{o}, \bm{y} }^{\bm{\varepsilon}} \!\geq\! \bm{0}, \bm{U}_{\rm c}^T \bm{\pi}_{\bm{o}, \bm{y} }^{\bm{\varepsilon}} \!=\! \bm{a} ~ \forall \bm{y}_{\bm{o}}^{\bm{\varepsilon}} \!\in\! \Theta
    \end{align}
\end{subequations}
where $\bm{\pi}_{\bm{o}, \bm{y} }^{\bm{\varepsilon}}$ are corresponding decision variables of $\bm{\pi}_{\bm{o} }^{\bm{\varepsilon}}$ for a particular assignment of $\bm{y}_{\bm{o}}^{\bm{\varepsilon}}$. Given that $\bm{o}$ are binary variables, (\ref{eq-3-6-function-A-equivalent-dual-emu:1}) can be linearized by applying the Big-M method as
\begin{subequations}\label{eq-3-6-A-bigM}
    \begin{align}
        &\!\!\! \varpi_{\bm{o}}^{\bm{\varepsilon}} \!\leq\! (\bm{t} \!-\!\bm{P} \bm{x} \!-\! \bm{Q} \bm{x}_{\bm{\varepsilon}} \!-\! \bm{U}_{\rm b} {\bm{y}}_{\bm{o}}^{\bm{\varepsilon}})^T \!\bm{\pi}_{\bm{o}, \bm{y} }^{\bm{\varepsilon}} \!-\!\! \bm{1}^T\! (\bm{O}^T  \!\!\circ\! \bm{\Pi}) \bm{1} \label{eq-3-6-A-bigM:1} \\
        &\!\!\!\!\! -\!\! M \bm{o}^T \bm{1}_{n_{\rm t}}  \leq_{[\bm{O}^T]} \bm{\Pi} \leq_{[\bm{O}^T]}  M \bm{o}^T \bm{1}_{n_{\rm t}} \label{eq-3-6-A-bigM:2} \\
        &\!\!\!\!\! -\!\! M \!(\!\bm{1}_N \!\!-\! \bm{o})^T \! \bm{1}_{n_{\rm t}}  \!\!\leq_{[\bm{O}^T]}\! \bm{\Pi} \!\!-\!\! \bm{1}_N^T \bm{\pi}_{\bm{o}, \bm{y} }^{\bm{\varepsilon}} \!\leq_{[\bm{O}^T]}\!  M (\!\bm{1}_N \!\!-\! \bm{o})^T \!\bm{1}_{n_{\rm t}} \!\!\! \label{eq-3-6-A-bigM:3}
    \end{align}
\end{subequations}
for each $\bm{y}_{\bm{o}}^{\bm{\varepsilon}} \!\in\! \Theta$. Here $\bm{\Pi} \!\in\! \mathbb{R}^{N \times n_{\rm t}}$ is the matrix of auxiliary variables, which is with the same sparsity as $\bm{O}^T$, i.e., only elements of  $\bm{\Pi}$ corresponding to non-zero elements of $\bm{O}^T$ are variables; and ``$\leq_{[\bm{O}^T]}$'' indicates that only inequalities corresponding to non-zero elements of $\bm{O}^T$ are valid. Then, (\ref{eq-3-6-function-A-equivalent-dual-emu}) can be solved as a mixed-integer linear program.

\begin{remark}
    The Big-M method is heavily used in the optimization model and solution approach, where improper values of Big-M parameters can potentially cut off legitimate solutions or cause numerical and convergence issues \cite{4-1322, 4-543}. In this work, we adopt an ad hoc Big-M parameter setting based on our experiences on solving OTS problems. Specifically, for the first kind of constraints which are those associated with physical requirements of transmission systems, including (\ref{eq-3-6-lpac:1})-(\ref{eq-3-6-lpac:6}), (\ref{eq-3-6-power-balance:1}), (\ref{eq-3-6-power-balance:2}), (\ref{eq-3-6-operation-1:5}), (\ref{eq-3-6-connected-1}), (\ref{eq-3-6-connected-zo})-(\ref{eq-3-6-phi-3-4-1}), (\ref{eq-3-6-connected-5:1}), and (\ref{eq-3-6-connected-5:2}), the Big-M values are all set as $10^7$; for the second kind of constraints which are associated with dual variables, i.e., (\ref{eq-3-6-A-bigM:2}), (\ref{eq-3-6-A-bigM:3}), and the second constraint in (\ref{eq-3-6-KKT-MI}), the Big-M values are all set as $10^8$; for the third kind of constraints which already contain other Big-M parameters, i.e., the first constraint in (\ref{eq-3-6-KKT-MI}), the Big-M values are set as 5 times of those already contained. In addition, taking CPLEX for example, the integrality tolerance is set to $10^{-15}$, and numerical precision emphasis is turned on to avoid potential numerical issues maximally. These settings lead to satisfactory performance in the case study. 
    It should be noted that this ad hoc setting can be improved by deriving tighter constraint-specific Big-M parameter values based on the knowledge of the first and second kinds of constraints above, and using the equivalent formulation with special ordered sets of type 1 for the third kind of constraints \cite{4-1273}. 
\end{remark}

\begin{remark}
    In practical applications, system operators can directly use the optimal solution of the TSDR-OTS model for the first-stage scheduling based on the forecast of VRE; then given realized VRE uncertainty, problem (\ref{eq-3-6-stage2-final}) is solved to yield the solution for the second-stage corrective control; when a contingency occurs, system operators implement the third-stage corrective control using the solution obtained by solving problem (\ref{eq-3-6-reform-ccg}). 
    Moreover, the solution of each stage can be AC-infeasible due to the power flow approximation, while common techniques are available to restore AC feasibility before applying the solution to real AC systems \cite[Ch. 6]{4-490}. 
\end{remark}

\section{Case Study}

The TSDR-OTS approach developed in this paper is implemented in Python using a Linux 64-Bit server with 2 Intel(R) Xeon(R) E5-2640 v4 @ 2.40GHz CPUs (a total of 40 threads provided) and 125GB RAM, where CPLEX 20.10 is used as solvers. Detailed description and data of all test systems used can be found in \cite{4-971-1}. Unless otherwise specified, lines 1 to 6, 24, and 25 of Algorithm \ref{alg-3-6-1} are executed using all threads of the server, and lines 7 to 23 are divided into 10 groups based on $\bm{\varepsilon}$ and each group is executed using 4 threads. 
The CPLEX optimality gap is set to 0.5\%. For Algorithm \ref{alg-3-6-1}, $\epsilon_1 = \epsilon_3 = 1\%$, $n_{\rm o} = 5$, $n_{\max} = 20$; $\epsilon_2 = 0.001 $ and $\epsilon_2 = 1 $ for the following simple IEEE network and real-scale networks, respectively. 
Since reliability data is not provided in the original data of the test systems, the parameters of contingency uncertainties are set based on the typical values of component failure probability \cite{4-1356}. Specifically, $(\bm{o}_{\min}, \bm{o}_{\max})$ is set as (0.0085, 0.0115) for generators, (0.004, 0.006) for transformers, and (0.00075, 0.00125) for transmission lines. Additionally, for each test system, we assume that each entry of $\bm{\varepsilon}$ follows the normal distribution with 0 mean and 15\% variance. To obtain $\Xi$, we first generate 200 and 5000 scenarios of $\bm{\varepsilon}$ by Monte Carlo simulations for the simple IEEE network and each real-scale network, respectively; then scenario reduction techniques in \cite{4-1521} are used to preserve 10 and 40 scenarios, respectively. Note that the number of preserved scenarios can be increased for higher solution accuracy given more computing units in practical applications.




\subsection{Case Study on the Simple IEEE Network}

The modified IEEE 9-bus system is first used to investigate the profit from line switching with the proposed TSDR-OTS model. 
This system is obtained from the original IEEE 9-bus system by adding 5 lines, and replacing two conventional generators with wind turbines. 
The dimensions of $\bm{\varepsilon}$ and $\bm{o}$ are 2 and 17, respectively. 
We consider 500 different scenarios of $\dot{\bm{p}}_{\rm d}$ and $\dot{\bm{p}}_{\rm v}^{\max}$; 
and different involvement of line switching, i.e., line switching at the first and third stages, line switching at the first stage only, line switching at the third stage only, and no line switching, which are denoted by S1 to S4, respectively. For comparison, we also consider different presence of uncertainties, i.e., both uncertainties, the uncertainty in forecast errors only, the uncertainty in contingencies only, and no uncertainties, which are denoted by U1 to U4, respectively. 

\begin{figure}[t]
	\centering
	\includegraphics[scale=0.8]{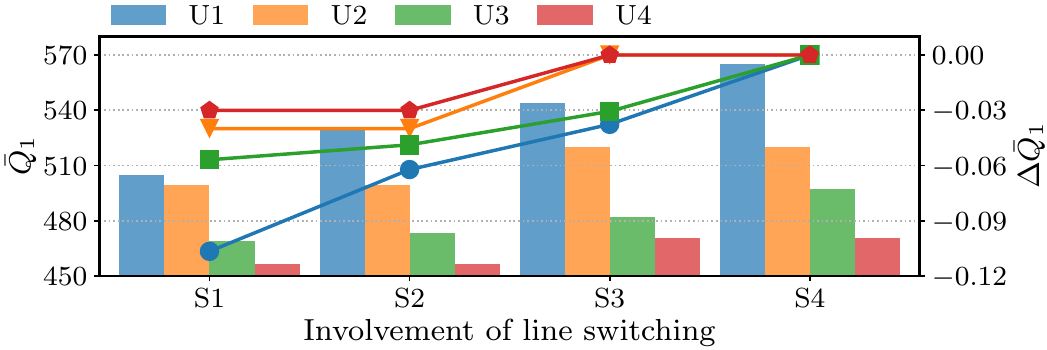}  
	\caption{The bar graph of $\bar{Q}_1$ and line graph of $\Delta \bar{Q}_1$ with various involvement of line switching and presence of uncertainties, where $\bar{Q}_1$ denotes the mean of $Q_1$ of the 500 different scenarios, and $\Delta \bar{Q}_1$ denotes the relative change from $\bar{Q}_1$ with no line switching and uncertainties to $\bar{Q}_1$ under certain case.}
	\label{fig-3-6-r2}
    \vspace*{-8pt}
\end{figure}

The TSDR-OTS model and its variations with the different involvement of line switching and presence of uncertainties are solved for the 500 scenarios, with the results given by Fig. \ref{fig-3-6-r2}. It is found that for any presence of uncertainties, $\bar{Q}_1$ increases overall with the involvement of line switching from S1 to S4, which indicates that more involvement of line switching can reduce the total operational cost and the reduction caused by the first-stage line switching is more significant than that caused by the third-stage one. The curve of $\Delta \bar{Q}_1$ moves up overall with the presence of uncertainties from U1 to U4, which implies that the profit from line switching is enlarged with more kinds of uncertainties involved. More importantly, this profit is maximized under U1 and S1, which is exactly the case corresponding to the proposed TSDR-OTS model.

\begin{figure}[t]
	\centering
 	\includegraphics[scale=0.8]{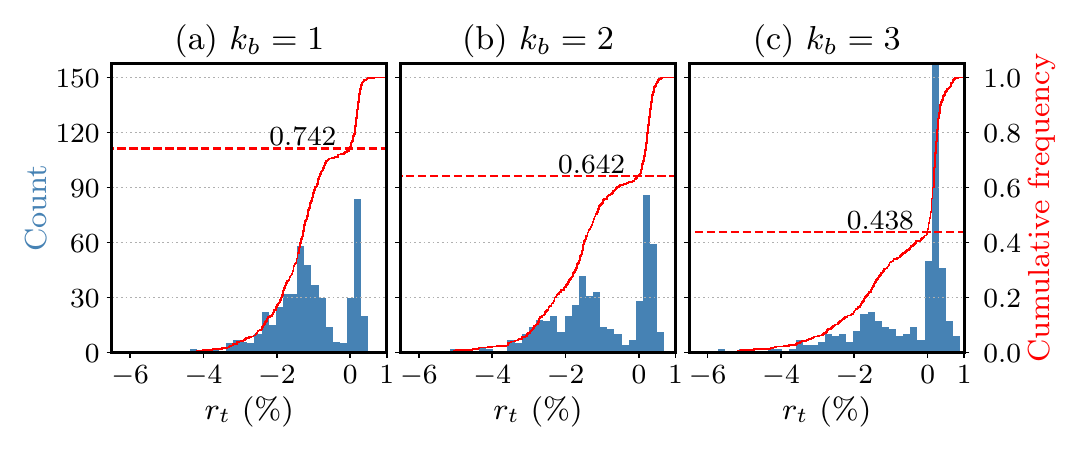}  
	\caption{Histograms of $r_t$ where $r_t = (t_1 - t_3)/{t_3}$ with $t_1$ and $t_3$ being solution time of M1 and M2, respectively.}
	\label{fig-3-6-r1}
\end{figure}

We further study the role of the first-stage NC constraints (\ref{eq-3-6-connected-1}). 
Consider the TSDR-OTS model and its variations with different consideration of NC, i.e., M1 to M4, referring to the TSDR-OTS models with full NC constraints, the first-stage NC constraints only, the third-stage NC constraints only and no NC constraints, respectively. 
Assume that $k_{\max} \!\!=\!\! 3$. Each model with $k_{\rm b} \!\!\in\!\! \{1,2,3\}$ is solved for the 500 scenarios. 
Table \ref{table-3-6-1} gives the NC ratios of the optimal topologies under the normal operation state and post-contingency state. 
It is found that for M1, $\zeta_{k} \!=\! 1$ for any $k \!\leq\! k_b$, which indicates that NC of the optimal topology under the normal operation state and after any $N\!-\!k_b$ contingency is always ensured. For M3 where (\ref{eq-3-6-connected-1}) is dropped, all its NC ratios equal to that of M1. In addition, for M2 and M4, the optimal topology under post-contingency state can be unconnected for both M2 and M4, while for M2, NC of the optimal topology under the normal operation state is always ensured. 
Therefore, the first-stage NC constraints are necessary when the third-stage ones are not contained, but are redundant in the TSDR-OTS model. 
Nevertheless, we find that these redundant constraints can improve computational efficiency overall. Fig. \ref{fig-3-6-r1} compares solution time of M1 and M3. When $k_b=1$, solution time of M1 is smaller than that of M3 for more than 70\% scenarios. With the increase of $k_b$, this proportion decreases while solution time of M1 is still smaller than that of M3 overall since the reduction of solution time caused by adding the first-stage NC constraints is generally more significant than the increase. The reason behind this counter-intuitive statistical finding can be explained by the structure of Algorithm \ref{alg-3-6-1}. With the third-stage NC constraints already contained, the first-stage ones are redundant for the entire TSDR-OTS model, but may not for MP $Q_1(\cdot)$ solved iteratively in Algorithm \ref{alg-3-6-1}. Hence, adding the first-stage NC constraints can speed up the convergence of Algorithm \ref{alg-3-6-1} for some scenarios.

\begin{table}[t]
	\centering
    \caption{Result statistics of network connectedness}
    \setlength{\tabcolsep}{0pt} 

    \setlength{\aboverulesep}{0pt}
    \setlength{\belowrulesep}{0pt}
    \setlength{\extrarowheight}{.1ex}
     \small{
    \begin{tabular*}{\hsize}{@{}p{0.77cm}@{}>{\columncolor{red!0}}P{0.62cm}@{}>{\columncolor{red!0}}P{0.6cm}@{}>{\columncolor{red!0}}P{0.6cm}@{}>{\columncolor{red!0}}P{0.6cm}@{}>{\columncolor{blue!0}}P{0.6cm}@{}>{\columncolor{blue!0}}P{0.6cm}@{}>{\columncolor{blue!0}}P{0.6cm}@{}>{\columncolor{blue!0}}P{0.6cm}@{}>{\columncolor{green!0}}P{0.6cm}@{}>{\columncolor{green!0}}P{0.6cm}@{}>{\columncolor{green!0}}P{0.6cm}@{}>{\columncolor{green!0}}P{0.6cm}}\toprule
        & \multicolumn{4}{c}{$k_{\rm b} = 1$ } & \multicolumn{4}{c}{$k_{\rm b} = 2$} & \multicolumn{4}{c}{$k_{\rm b} = 3$}  
        \\\cmidrule(l){2-5}\cmidrule(l){6-9}\cmidrule(lr){10-13} 
               & $\zeta_{0}$ & $\zeta_1$ & $\zeta_2$ & $\zeta_3$  & $\zeta_{0}$ & $\zeta_1$ & $\zeta_2$ & $\zeta_3$  & $\zeta_{0}$ & $\zeta_1$ & $\zeta_2$ & $\zeta_3$  \\ \midrule
               M1 & 1 & 1 & 0.96 & 0.88 & 1 & 1 & 1 & 0.99 & 1 & 1 & 1 & 1 \\
               M2 & 1 & 0.97 & 0.90 & 0.83 & -- & -- & -- & -- & -- & -- & -- & -- \\
               M3 & 1 & 1 & 0.96 & 0.88 & 1 & 1 & 1 & 0.99 & 1 & 1 & 1 & 1 \\
               M4 & 0.96 & 0.93 & 0.86  & 0.80 & -- & -- & -- & -- & -- & -- & -- & --\\  \bottomrule
    \multicolumn{13}{@{}p{1\columnwidth}@{}}{
        \footnotesize{\textit{Note}: 
        $\zeta_0$ is the ratio of the optimal topology that is connected under the normal operation state, and $\zeta_k$ with $k \!\in\! \{1,2,3\}$ is the ratio of pairs of the optimal topology and $N\!\!-\!\!k$ contingency for which the number of connected components of the corresponding post-contingency topology equals to that of the post-contingency topology induced from the same contingency and the topology with all lines closed. Their computation only considers contingencies consisting of branch failures entirely among all possible $N-k$ contingencies.
      }}
    \end{tabular*} 
     }
    \label{table-3-6-1}  
\end{table}

\subsection{Case Study on the Real-Scale Networks}


Practical effectiveness of the proposed TSDR-OTS model and efficiency of the solution approach are further tested with the 50Hertz control area of the German transmission networks (DE-50Hz for short) and the transmission network of Liaoning province in China (LN for short). 
For the DE-50Hz (LN) network, the VRE penetration rate is 38.49\% (19.86\%); the first-stage model contains 292 (312) continuous decision variables and 120 (200) binary decision variables, the second-stage model contains 438 (468) continuous decision variables, the third-stage model contains 505 (766) continuous decision variables and 240 (400) binary decision variables; the maximal number of fault components $k_{\max}$ is set as 4 (3), and $k_{\rm b}$ is set as 2 (2); the dimensions of $\bm{\varepsilon }$ and $\bm{o}$ are 57 (38) and 258 (561), respectively.  
It is noted that for the LN network, we assume 200 out of 954 transmission lines participate in OTS and the third-stage corrective control, and constraints $\bm{1}^T \dot{\bm{z}} \!\!\geq\!\! 992$ and $\bm{1}^T\! \bm{z}_+ \!\!+\!\! \bm{1}^T\! \bm{z}_- \!\!\leq\!\! 10$ are added to the TSRD-OTS model.

\subsubsection{Effectiveness of the TSDR-OTS approach}

\begin{figure}[t]
    \centering
    \begin{subfigure}[b]{0.85\linewidth}
        \includegraphics[width=\linewidth]{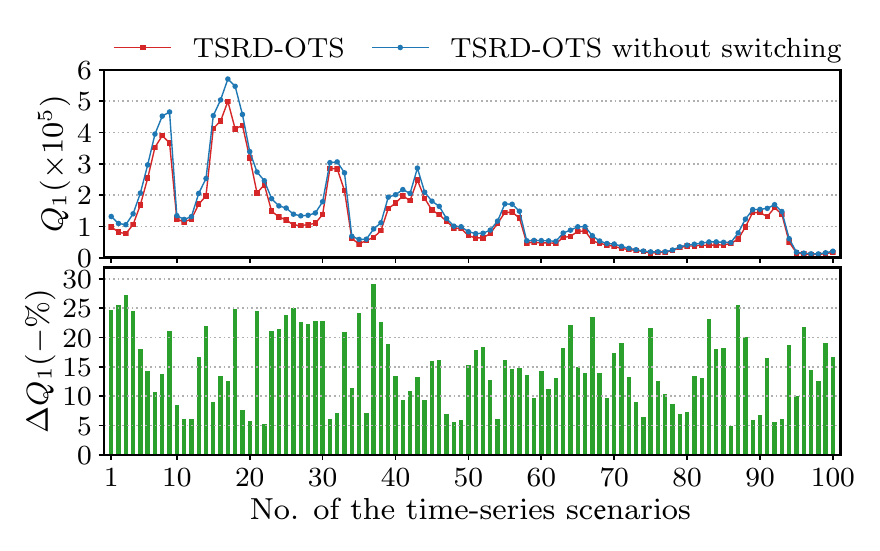} 
        \caption{\footnotesize{DE-50Hz}}
        \label{fig-4-1-dnr-r21:evaluation}
    \end{subfigure} 
    \begin{subfigure}[b]{0.85\linewidth}
        \includegraphics[width=\linewidth]{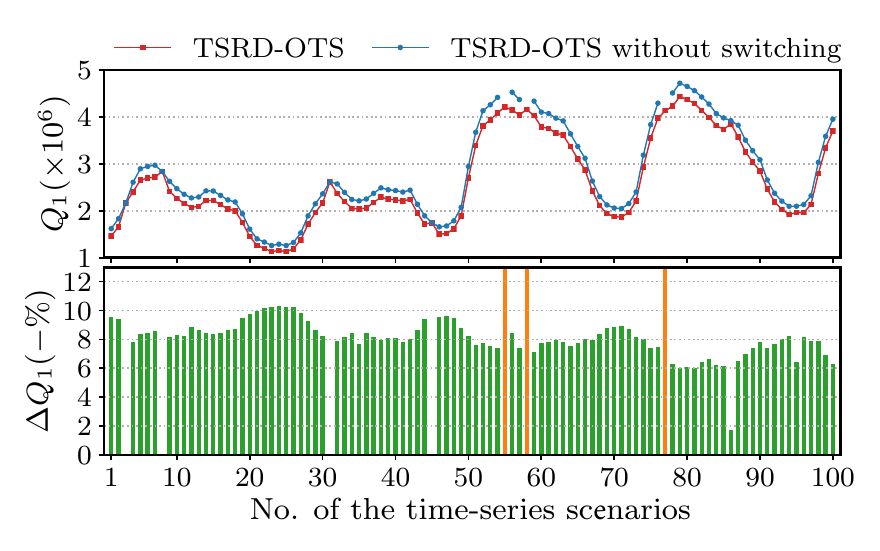} 
        \caption{\footnotesize{LN}}
        \label{fig-4-1-dnr-r21:time}
    \end{subfigure}
    \caption{Comparison of $Q_1$ given by the TSRD-OTS model and that without line switching for 100 time-series scenarios (the upper 2 figures), and the corresponding relative change from the latter $Q_1$ to the former one (the lower 2 figures), i.e., $\Delta Q_1$. 
    }\label{fig-3-6-r-3-4} 
\end{figure}

For each system, the proposed TSRD-OTS model and its variation where line switching is not allowed are both solved for 100 time-series scenarios with a 1-hour interval, with the results given by Fig. \ref{fig-3-6-r-3-4}. 
In Fig. \ref{fig-3-6-r-3-4} (b), the missing points in the blue curve indicate that values of $Q_1$ are abnormally large or the solver finds the model infeasible in line 3 of Algorithm \ref{alg-3-6-1}. The orange bars mean invalid values of $\Delta Q_1$. 
Moreover, 
Table \ref{table-3-6-add} reports the aggregate statistics for the obtained optimal solutions and computational performance, where $c_{\rm se}$ and $c_{\rm se}'$ indicate the conservativeness of the solution to ensure the $N\!-\!k$ security, and all statistics are for the TSRD-OTS model except for $c_{\rm se}'$ that is for the above model variation.

According to Fig. \ref{fig-3-6-r-3-4}, we can find that for the DE-50Hz network, the proposed TSRD-OTS effectively reduces $Q_1$ for all the scenarios compared with the case where line switching is not allowed. The maximal rate of decrease is close to 30\%, while the minimum is still about 5\%. For the LN network, the proposed TSRD-OTS reduce $Q_1$ by 6\% to 10\% for the majority of the scenarios. In particular, for some scenarios, i.e., the 55th, 58th, and 77th ones, no feasible solution exists when the network topology is fixed, while the proposed TSRD-OTS ensures the feasibility by leveraging the topological flexibility. 
Moreover, the value of $c_{\rm se}$ is smaller than $c_{\rm se}'$ for all cases in Table \ref{table-3-6-add}, indicating that the first-stage line switching reduces the conservativeness of the solution to ensure the $N\!-\!k$ security. 
For the solution time, the average even for the DE-50Hz network is more than 60 minutes, which was expected given the complexity of the TSRD-OTS model. 
Nevertheless, the computational suitability of our approach for practical applications can be easily ensured by using more computing units to execute the DW procedure of each $\bm{\varepsilon} \in \Xi$ in parallel.

\begin{table}[t]
	\centering
    \caption{Statistics for the solutions and computations.}
    \setlength{\tabcolsep}{0pt} 

    \setlength{\aboverulesep}{0pt}
    \setlength{\belowrulesep}{0pt}
    \setlength{\extrarowheight}{.2ex}
     \small{
    \begin{tabular*}{\hsize}{@{}p{0.8cm}@{}p{1.6cm}@{}p{1cm}@{}p{1.2cm}@{}p{1.2cm}@{}p{1.1cm}@{}p{1cm}@{}p{1cm}}\toprule
               &  & $n_{\rm sl}$  & $c_{\rm se}$ & $c_{\rm se}'$ &  $t_{\rm st}$ & $n_{\rm in}$  & $n_{\rm out}$   \\ \midrule
               \parbox[t]{1mm}{\multirow{2}{*}{~$\!$\rotatebox[origin=c]{90}{\footnotesize{  \makecell{DE- \\ 50Hz} }}}} 
                & Average     & 9.41   & 12963  & 17046   & 73.83  & 7.34  &  9.72   \\ 
                & Worst-case  & 16      & 22626  & 23058  & 95.42  &  12 & 19     \\ \midrule
               \parbox[t]{1mm}{\multirow{2}{*}{~$\!$\rotatebox[origin=c]{90}{LN}}}
                & Average    &  12.89   & 185231  & 235703  & 184.34  & 8.53  & 8.03 \\ 
                & Worst-case &  23      & 311411  & 344659  & 242.65  & 14  & 23 \\ \bottomrule
            \multicolumn{8}{@{}p{1\columnwidth}@{}}{
                \footnotesize{\textit{Note}: 
                $n_{\rm sl}$ denotes the the number of switched lines of the first stage, $t_{\rm st}$ the total solution time in minutes; $n_{\rm in}$ and $n_{\rm out}$ the iteration number of the inner and outer CCG loop, respectively; 
                $c_{\rm se}$ and $c_{\rm se}'$ the difference between $Q_1$ with the third-stage objective ignored of the optimal solution of the TSRD-OTS model and that of the TSRD-OTS model with the third stage excluded.
            }}
    \end{tabular*} 
     }
    \label{table-3-6-add}  
\end{table}

\subsubsection{AC feasbility of the optimal solutions}
Considering the utilization of the approximated power flow model, we further study the AC feasibility of the obtained optimal solutions and the impact of potential AC restoration on the evaluated objective function values. 
Specifically, the AC feasibility is evaluated based on all the 100 time-series scenarios and all samples of $\bm{\varepsilon}$ for the first- and second-stage decisions, and 20 time-series scenarios, 10 samples of $\bm{\varepsilon}$ and 500 randomly selected contingencies for the third-stage decisions. For AC-infeasible decisions, we use the sequential quadratic programming algorithm to restore their AC feasibility. 
Table \ref{table-3-6-4} reports the aggregate statistics for the decisions of each stage. It presents the ratio of AC-infeasible decisions ($r_{\rm{inf}}$), the average ($\mu$) and worst-case ($\max$) values of the absolute error ($Q_{\Delta}$) between the objective value of the AC-infeasible decision and the restored AC-feasible one, the average relative error ($\mu(\delta(Q_{\Delta}))$) and the relative error of the value selected by the max operator ($\delta({\max}(Q_{\Delta})) $).

\begin{table}[t]
	\centering
    \caption{Statistics for AC feasibility and its restoration.}
    \setlength{\tabcolsep}{0pt} 

    \setlength{\aboverulesep}{0pt}
    \setlength{\belowrulesep}{0pt}
    \setlength{\extrarowheight}{.2ex}
     \small{
    \begin{tabular*}{\hsize}{@{}p{0.5cm}@{}p{1.2cm}@{}p{1cm}@{}p{1.1cm}@{}p{1.57cm}@{}p{1.57cm}@{}p{1.37cm}}\toprule
               &  Stage   & $r_{\rm{inf}}$ & $\mu(Q_{\Delta})$ &  $\mu(\delta(Q_{\Delta}))$ & $\max(Q_{\Delta})$  & $\delta({\max}(Q_{\Delta})) $   \\ \midrule
               \parbox[t]{1mm}{\multirow{3}{*}{~$\!$\rotatebox[origin=c]{90}{\footnotesize{DE-50Hz}}}} 
                &  First   & 0.950  & 565  & 0.0071 & 2168  &   0.033    \\
                &  Second  & 0.852  & 125  & 0.0034 & 363   &   0.027     \\
                &  Third   & 0.751  & 153 & 0.018  & 805  &   0.093     \\ \midrule
               \parbox[t]{1mm}{\multirow{3}{*}{~$\!$\rotatebox[origin=c]{90}{LN}}}
                &  First   & 0.930  & 7693   &   0.0042 & 216504  & 0.091 \\
                &  Second  & 0.947  & 4106   &   0.0063 & 18837   & 0.046  \\
                &  Third   & 0.811  & 2150  &   0.015  & 13650  & 0.13  \\ \bottomrule
    \end{tabular*} 
     }
    \label{table-3-6-4}  
\end{table}



Table \ref{table-3-6-4} indicates that most of the obtained decisions are feasible to the AC system, especially for the first stage. For the AC-infeasible decisions of each stage, their AC feasibility can be restored with slight changes of the objective function values. Moreover, the benefit from line switching evaluated via the approximated model is still valid overall after the AC feasibility restoration, according to the values of $\mu(\delta(Q_{\Delta}))$.

\subsubsection{Comparison with alternative solution algorithms}

Regarding computational efficiency, we compare solution time of using Algorithm \ref{alg-3-6-1} and some alternatives for the DE-50Hz network's first 50 scenarios. Meanwhile, impacts of constraints (\ref{eq-3-6-connected-1}) and (\ref{eq-3-6-symmetry}) on solution time are also analysed. Consider Algorithm \ref{alg-3-6-1}, the nested CCG algorithm in \cite{4-1265} , Algorithm \ref{alg-3-6-1} with the inner loop's master problem derived using the KKT condition, Algorithm \ref{alg-3-6-1} with all lines executed using all threads of the server, Algorithm \ref{alg-3-6-1} solving the TSDR-OTS model without (\ref{eq-3-6-connected-1}), and Algorithm \ref{alg-3-6-1} solving the TSDR-OTS model without (\ref{eq-3-6-symmetry}), which are called A1 to A6, respectively. Table \ref{table-3-6-2} gives the results of comparison. The results of A2 and A3 indicate that the DW procedure and reformulating $\mathcal{A}$ based on strong duality both significantly reduce solution time of using Algorithm \ref{alg-3-6-1}. The result of A5 indicates that constraints (\ref{eq-3-6-connected-1}) can improve computational efficiency overall, although they are redundant. The result of A6 demonstrates the potential reduction of solution time by adding constraint (\ref{eq-3-6-symmetry}). In particular, the most significant improvement of computational efficiency is achieved by parallel computing for the DW procedure. Overall, Algorithm 1 outperforms all the alternatives regarding computational efficiency.

\begin{table}[t]
    \caption{Comparison of solution time using A1 to A6.}
    \centering
    \setlength{\extrarowheight}{-.0ex}
    \footnotesize{
    \begin{tabularx}{1.0\columnwidth}{@{} p{1.4cm}<{\centering}p{1.4cm}<{\centering}p{1.4cm}<{\centering}p{1.4cm}<{\centering}p{1.4cm}<{\centering} @{}}
      \toprule 
      A2 & A3 & A4 & A5 & A6 \\[-0.6ex]
      \midrule
      \begin{minipage}{0.95\textwidth} \includegraphics[width=0.095\textwidth]{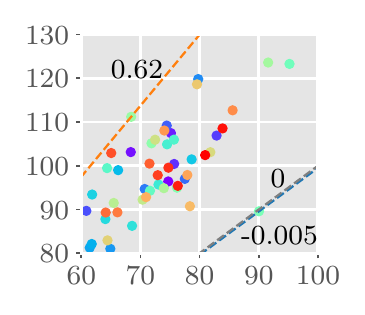} \end{minipage} &
      \begin{minipage}{0.95\textwidth} \includegraphics[width=0.095\textwidth]{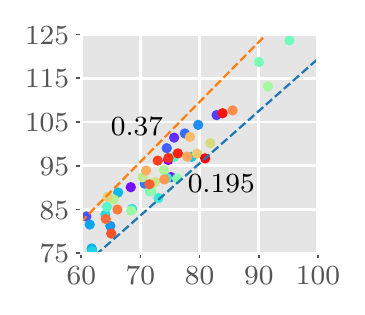} \end{minipage} &
      \begin{minipage}{0.95\textwidth} \includegraphics[width=0.095\textwidth]{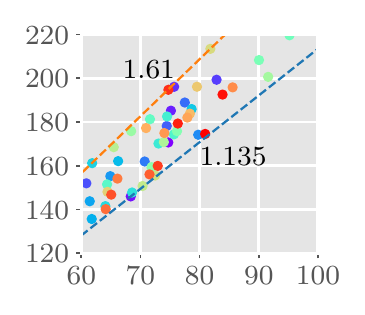} \end{minipage} &
      \begin{minipage}{0.95\textwidth} \includegraphics[width=0.095\textwidth]{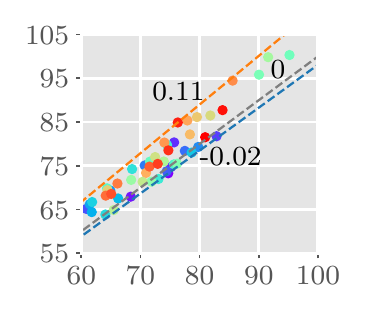} \end{minipage} &
      \begin{minipage}{0.95\textwidth} \includegraphics[width=0.095\textwidth]{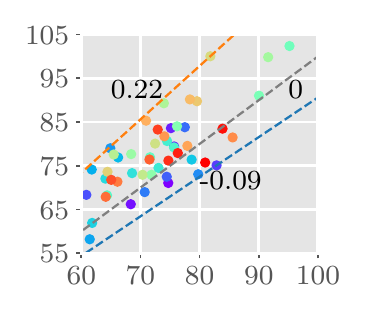} \end{minipage} \\
      \midrule
      35.73\% & 27.73\% & 134.42\% & 4.22\% & 6.69\% \\[-0.6ex]
      \bottomrule
      \multicolumn{5}{@{}p{1\columnwidth}@{}}{
        \footnotesize{\textit{Note}: The second row is the $\tau_1$-$\tau_i$ scatter plot for the first 50 scenarios of the DE-50Hz network where $\tau_i$ is solution time (in minutes) using A$i$ given as the column name. The third row is the mean of $\Delta \tau_i = (\tau_i - \tau_1)/\tau_1$. The number next to a dash line is the value of $\Delta \tau_i$ for points on the line.
        }}
      \end{tabularx}}
	\label{table-3-6-2}
    \vspace{-5pt}
\end{table}

\subsubsection{Comparison with two-stage stochastic OTS}

We compare the TSDR-OTS approach to the two-stage stochastic OTS (TSS-OTS) established under the framework in \cite{4-1139} and \cite{4-1142}. 
Under this framework, the first stage is the same as that in the TSDR-OTS model; the second stage is the corrective control reacting to probabilistic scenarios representing simultaneous occurrence of VRE forecast errors and contingencies. 
The second-stage corrective control in the TSS-OTS model is set as the third-stage corrective control in the TSDR-OTS model since otherwise the TSS-OTS model is unobtainable. 
Only 10 $N\!-\!1$ contingencies with the probabilities of component failures being $(\bm{o}_{\min} \!+\! \bm{o}_{\max})/2$ are considered in the TSS-OTS model for computational tractability and the same contingencies plus the normal state are used as the support $\Omega$ in the TSDR-OTS model for comparability. Meanwhile, the values of $\bm{o}_{\min}$ and $\bm{o}_{\max}$ are increased to compensate for the impact of reducing $\Omega$ (see \cite{4-971-1} for detailed data). 
Assume that $O^*: \mathbb{E}_{O}(\bm{1} - \bm{o}) = (\bm{o}_{\min} + \bm{o}_{\max})/2$ is the true PD of $\bm{o}$. Considering that the VRE uncertainty and contingency happen separately in practical, the performance of each OTS model is evaluated as follows: solve the OTS model to obtain the first-stage scheduling solution with the associate objective value being $c_1$; with this scheduling solution, for each $\bm{\varepsilon} \in \Xi$, solve the second-stage corrective control problem to obtain the associated regulation cost denoted by $c_{2}^{\bm{\varepsilon}}$; following the second-stage corrective control induced by $\bm{\varepsilon}$, for each $\bm{o} \in \Omega$, solve the third-stage corrective control problem to obtain the associated corrective control cost $c_{3}^{\bm{\varepsilon}, \bm{o}}$; compute $Q_1^* \!=\! c_1 \!+\! \sum_{\bm{\varepsilon} \in \Xi} P(\bm{\varepsilon}) [ c_{2}^{\bm{\varepsilon}} \!+\! \sum_{\bm{o} \in \Omega} O^*(\bm{\varepsilon}) c_{3}^{\bm{\varepsilon}, \bm{o}} ]$, which is the total operational cost with the true PD of $\bm{o}$ and practical sequence of uncertainty occurrence.

\begin{table}[t]
	\centering
    \caption{Values of $Q_1^*$ of different OTS models.}
    \setlength{\tabcolsep}{0pt} 

    \setlength{\aboverulesep}{0pt}
    \setlength{\belowrulesep}{0pt}
    \setlength{\extrarowheight}{.2ex}
     \small{
    \begin{tabular*}{\hsize}{@{}p{0.48cm}@{}p{1.7cm}@{}p{1.37cm}@{}p{1.37cm}@{}p{1.37cm}@{}p{1.37cm}@{}p{1.15cm}}\toprule
        \parbox[t]{1mm}{\multirow{3}{*}{~$\!$\rotatebox[origin=c]{90}{\footnotesize{DE-50Hz}}}} 
               &  Model   & No.7 & No.8 & No.22  & No.29 & No.47 \\ \cmidrule{2-7}
               &  TSDR-OTS  & 382314  & 348972  & 138533   & 131472   & 92456  \\
               &  TSS-OTS   & 397881  & 357411  & 144780   & 138409   & 96980   \\ \toprule
               \parbox[t]{1mm}{\multirow{3}{*}{~$\!$\rotatebox[origin=c]{90}{LN}}}
                &  Model   & No.4 & No.25 & No.42  & No.77 & No.81 \\ \cmidrule{2-7}
                &  TSDR-OTS & 2600620 & 1129983 & 1817773  & 4023560  & 4014464  \\
                &  TSS-OTS  & 2663538 & 1162542 & 1844586  & 4181332  & 4080097  \\ \bottomrule
    \end{tabular*} 
     }
    \label{table-3-6-3}  
\end{table}

Table \ref{table-3-6-3} compares the values of $Q_1^*$ associated with different OTS models, for some time-series scenarios of the DE-50Hz and LN networks. The true PD of $\bm{o}$ is employed by the TSS-OTS model but is not necessarily the worst-case one associated with the optimal first-stage scheduling solution yield by the TSDR-OTS model. Thus, $Q_1^*$ of the TSS-OTS model should not be larger than that of the TSDR-OTS model for the same time-series scenario. However, for each time-series scenario in Table \ref{table-3-6-3}, the opposite of this inference is observed. This indicates that with the practical sequence of uncertainty occurrence and separate corrective controls, the TSS-OTS approach can lead to inferior first-stage scheduling solutions compared to the TSDR-OTS approach.

\subsubsection{Comparison with alternative modelling of contingency uncertainty}

Finally, we study the superiorities of the DRO technique used for tackling contingency uncertainty over the alternatives including the RO and SO techniques. 
With RO, only the support of $\bm{o}$ is utilized to model its uncertainty and the term $\sup_{ O \in \mathcal{O} }   \mathbb{E}_{ O } [ Q_3( \bm{x}_{\bm{\varepsilon}},\! \bm{x},\! \bm{o} ) ]$ in (\ref{eq-3-6-final-stage-2:1}) is replaced by $\sup_{ \bm{o} \in \Omega }  Q_3( \bm{x}_{\bm{\varepsilon}},\! \bm{x},\! \bm{o} )$. With SO, this term is replaced by $\mathbb{E}_{ O } [ Q_3( \bm{x}_{\bm{\varepsilon}},\! \bm{x},\! \bm{o})]$ with $O$ being the assumed exact PD of $\bm{o}$. All the other factors are equal except for these differences. It is assumed that $O^*$ is the exact PD of $\bm{o}$ in SO. For computational tractability, we only consider all the $N-k$ contingencies over 5 lines and 5 generators. 
Table \ref{table-3-6-5} reports 
the third-stage corrective control cost (TCC) under the worst-case contingency (WCC), 
the expected TCC under the worst-case PD of contingency uncertainty (WCD), 
out-of-sample expected TCC under a randomly simulated PD of contingency uncertainty from $\mathcal{O}$ (RSD); 
the total operational cost (TOC) with these different cases of TCC; 
and the computational time of each model. 
To obtain these values for the DRO, RO and SO models, we first solve the three models and their associated subproblems to obtain 
(i) the optimal first-stage decisions and the optimal second-stage decisions for each sample of $\bm{\varepsilon}$, of the three models; 
(ii) the worst-case contingency for each sample of $\bm{\varepsilon}$, from the RO model; 
(iii) the worst-case PD of contingency uncertainty for each sample of $\bm{\varepsilon}$, from the DRO model. 
Then for each sample of $\bm{\varepsilon}$, with the first- and second-stage decisions fixed, we evaluate different cases of TCC by solving problem (\ref{eq-3-6-reform}) for all possible $N-k$ contingencies. The time-series scenario 1 of the DE-50Hz network is used for the comparison here. All TCC values are averaged over all samples of $\bm{\varepsilon}$.

Table \ref{table-3-6-5} shows that the DRO technique yields not only lower third-stage corrective control cost but also lower total operational cost, both under the worst-case PD of contingency uncertainty and in the out-of-sample simulation. 
The RO model outperforms regarding both TCC and TOC for the worst-case contingency, which was expected, since RO optimizes the first- and second-stage decisions w.r.t. the worst-case contingency. Moreover, the solution time of each model demonstrates the superiority of the DRO model, particularly over the SO model, in terms of computational efficiency.

\begin{table}[t]
	\centering
    \caption{Comparison of TCC and TOC.}
    \setlength{\tabcolsep}{0pt} 

    \setlength{\aboverulesep}{0pt}
    \setlength{\belowrulesep}{0pt}
    \setlength{\extrarowheight}{.1ex}
     \small{
    \begin{tabular*}{\hsize}{@{}p{1cm}@{}p{1.1cm}@{}p{1.1cm}@{}p{1.1cm}@{}p{1.1cm}@{}p{1.1cm}@{}p{1.1cm}@{}p{1.1cm}}\toprule
        & \multicolumn{3}{c}{TCC } & \multicolumn{3}{c}{TOC } &  \\ \cmidrule(r){2-4}\cmidrule(r){5-7} 
      Model   & WCC   & WCD   &  RSD  &  WCC  & WCD   &  RSD  & Time (s)   \\ \midrule 
         SO   & 13917 & 6526 & 4869  & 88862 & 81471 & 79814 & 13490  \\                 
         RO   & 9452  & 7037 & 6533  & 83523 & 81108 & 80604 & 2462   \\                
         DRO  & 13034 & 4941 & 4271  & 88430 & 80337 & 79000 & 1142    \\ \bottomrule   
    \end{tabular*} 
     }
    \label{table-3-6-5}  
\end{table}

\section{Conclusion}

In this paper, we promoted filling the gap of the OTS problem of VRE penetrated power grids: how to simultaneously handle the VRE uncertainty and $N\!\!-\!k$ security in a proper and efficient way. 
The proposed TSDR-OTS model handles the uncertainties of VRE and contingencies based on SO and DRO respectively, given different availability of their PDs. In particular, stage-wise realization of the two kinds of uncertainties enables their corrective controls with different mechanisms to be considered separately and thus properly. With tractable reformulation derived, we developed an efficient solution approach whose effectiveness was demonstrated with different networks under various scenarios. Reduction of the total operational cost and recovery of solution feasibility by the proposed TSDR-OTS were also verified.

Future directions are twofold. 
Firstly, the proposed TSDR-OTS approach can be extended to OTS problems with unit commitment decisions and multi-period settings. However, the optimization model structure will be complicated with more stages involved, especially when the non-anticipativity of the corrective control decisions and the uncertainty in $\dot{\bm{p}}_{\rm v}^{\max}$ at the day-ahead stage are considered. Thus, challenges in computational efficiency need to be further tackled, for example by combining decomposition techniques leveraging the multi-period problem structure and affine decision rules. 
Secondly, since the third-stage post-contingency controls are not applied instantaneously in practical, the post-contingency system should be able to survive the time before the third-stage controls are applied. Thus the proposed approach can be improved by considering an additional stage in the formulation, with coupling constraints between the pre-contingency and every possible post-contingency stage. The resulting more complicated yet more practical formulation and its solution approach deserve further studies.

\ifCLASSOPTIONcaptionsoff
  \newpage
\fi

\bibliographystyle{IEEEtran}
\bibliography{References/4}

\begin{thebibliography}{10}
\providecommand{\url}[1]{#1}
\csname url@samestyle\endcsname
\providecommand{\newblock}{\relax}
\providecommand{\bibinfo}[2]{#2}
\providecommand{\BIBentrySTDinterwordspacing}{\spaceskip=0pt\relax}
\providecommand{\BIBentryALTinterwordstretchfactor}{4}
\providecommand{\BIBentryALTinterwordspacing}{\spaceskip=\fontdimen2\font plus
\BIBentryALTinterwordstretchfactor\fontdimen3\font minus \fontdimen4\font\relax}
\providecommand{\BIBforeignlanguage}[2]{{%
\expandafter\ifx\csname l@#1\endcsname\relax
\typeout{** WARNING: IEEEtran.bst: No hyphenation pattern has been}%
\typeout{** loaded for the language `#1'. Using the pattern for}%
\typeout{** the default language instead.}%
\else
\language=\csname l@#1\endcsname
\fi
#2}}
\providecommand{\BIBdecl}{\relax}
\BIBdecl

\bibitem{4-995-ea}
T.~Han, Y.~Song, and D.~J. Hill, ``Ensuring network connectedness in optimal transmission switching problems,'' \emph{{IEEE} Trans. Circuits Syst. {II}}, vol.~68, no.~7, pp. 2603--2607, Jul. 2021.

\bibitem{4-970}
J.~Li, F.~Liu, Z.~Li, C.~Shao, and X.~Liu, ``Grid-side flexibility of power systems in integrating large-scale renewable generations: A critical review on concepts, formulations and solution approaches,'' \emph{Renew. Sust. Energ. Rev.}, vol.~93, pp. 272--284, Oct. 2018.

\bibitem{4-62}
E.~B. Fisher, R.~P. Oneill, and M.~C. Ferris, ``Optimal transmission switching,'' \emph{IEEE Trans. Power Syst.}, vol.~23, no.~3, pp. 1346--1355, Aug. 2008.

\bibitem{4-670}
K.~W. Hedman, M.~C. Ferris, R.~P. O'Neill, E.~B. Fisher, and S.~S. Oren, ``Co-optimization of generation unit commitment and transmission switching with {N}-1 reliability,'' \emph{IEEE Trans. Power Syst.}, vol.~25, no.~2, pp. 1052--1063, May 2010.

\bibitem{4-63}
K.~W. Hedman, R.~P. Oneill, E.~B. Fisher, and S.~S. Oren, ``Optimal transmission switching with contingency analysis,'' \emph{IEEE Trans. Power Syst.}, vol.~24, no.~3, pp. 1577--1586, May 2009.

\bibitem{4-1140}
T.~Ding and C.~Zhao, ``Robust optimal transmission switching with the consideration of corrective actions for {N}-k contingencies,'' \emph{IET Gener. Transm. Distrib.}, vol.~10, no.~13, pp. 3288--3295, Oct. 2016.

\bibitem{4-1137}
F.~Qiu and J.~Wang, ``Chance-constrained transmission switching with guaranteed wind power utilization,'' \emph{IEEE Trans. Power Syst.}, vol.~30, no.~3, pp. 1270--1278, May 2015.

\bibitem{4-79}
P.~Dehghanian and M.~Kezunovic, ``Probabilistic decision making for the bulk power system optimal topology control,'' \emph{IEEE Trans. Smart Grid}, vol.~7, no.~4, pp. 2071--2081, Jul. 2016.

\bibitem{4-1223}
H.~Zhang, H.~Cheng, and S.~Zhang, ``Stochastic optimal transmission switching considering the correlated wind power,'' \emph{IET Gener. Transm. Distrib.}, vol.~13, no.~13, pp. 2664--2672, Jul. 2019.

\bibitem{4-1158}
M.~Alhazmi, P.~Dehghanian, S.~Wang, and B.~Shinde, ``Power grid optimal topology control considering correlations of system uncertainties,'' \emph{{IEEE} Trans. Ind. Appl.}, vol.~55, no.~6, pp. 5594--5604, Nov. 2019.

\bibitem{4-1139}
A.~Nikoobakht, J.~Aghaei, and M.~Mardaneh, ``Securing highly penetrated wind energy systems using linearized transmission switching mechanism,'' \emph{Applied Energy}, vol. 190, pp. 1207--1220, Mar. 2017.

\bibitem{4-1142}
A.~Nikoobakht, J.~Aghaei, M.~Mardaneh, T.~Niknam, and V.~Vahidinasab, ``Moving beyond the optimal transmission switching: stochastic linearised {SCUC} for the integration of wind power generation and equipment failures uncertainties,'' \emph{IET Gener. Transm. Distrib.}, vol.~12, no.~15, pp. 3780--3792, Aug. 2018.

\bibitem{4-1226}
A.~Nikoobakht, J.~Aghaei, M.~Lotfi, J.~P. Catalao, G.~J. Osorio, and M.~Shafie-khah, ``Flexible co-operation of {TCSC} and corrective topology control under wind uncertainty: An interval-based robust approach,'' in \emph{2019 {IEEE} Milan {PowerTech}}.\hskip 1em plus 0.5em minus 0.4em\relax {IEEE}, Jun. 2019.

\bibitem{4-1127}
H.~Huang, M.~Zhou, S.~Zhang, L.~Zhang, G.~Li, and Y.~Sun, ``Exploiting the operational flexibility of wind integrated hybrid {AC}/{DC} power systems,'' \emph{IEEE Trans. Power Syst.}, vol.~36, no.~1, pp. 818--826, Jan. 2021.

\bibitem{4-1248}
R.~Saavedra, A.~Street, and J.~M. Arroyo, ``Day-ahead contingency-constrained unit commitment with co-optimized post-contingency transmission switching,'' \emph{IEEE Trans. Power Syst.}, vol.~35, no.~6, pp. 4408--4420, Nov. 2020.

\bibitem{4-1195}
S.~Babaei, R.~Jiang, and C.~Zhao, ``Distributionally robust distribution network configuration under random contingency,'' \emph{IEEE Trans. Power Syst.}, vol.~35, no.~5, pp. 3332--3341, Sept. 2020.

\bibitem{4-1214}
C.~Zhao and R.~Jiang, ``Distributionally robust contingency-constrained unit commitment,'' \emph{IEEE Trans. Power Syst.}, vol.~33, no.~1, pp. 94--102, Jan. 2018.

\bibitem{4-538}
A.~S. Korad and K.~W. Hedman, ``Robust corrective topology control for system reliability,'' \emph{IEEE Trans. Power Syst.}, vol.~28, no.~4, pp. 4042--4051, Nov. 2013.

\bibitem{4-779}
X.~Li and K.~W. Hedman, ``Enhanced energy management system with corrective transmission switching strategy-{P}art {I}: {M}ethodology,'' \emph{IEEE Trans. Power Syst.}, vol.~34, no.~6, pp. 4490--4502, Nov. 2019.

\bibitem{4-1217}
X.~Lu, K.~W. Chan, S.~Xia, B.~Zhou, and X.~Luo, ``Security-constrained multiperiod economic dispatch with renewable energy utilizing distributionally robust optimization,'' \emph{IEEE Trans. Sustain. Energy}, vol.~10, no.~2, pp. 768--779, Apr. 2019.

\bibitem{4-1197}
P.~Li, Q.~Wu, M.~Yang, Z.~Li, and N.~Hatziargyriou, ``Distributed distributionally robust dispatch for integrated transmission-distribution systems,'' \emph{IEEE Trans. Power Syst.}, vol.~36, no.~2, pp. 1193--1205, Mar. 2021.

\bibitem{4-1228}
W.~Zheng, W.~Huang, D.~J. Hill, and Y.~Hou, ``An adaptive distributionally robust model for three-phase distribution network reconfiguration,'' \emph{IEEE Trans. Smart Grid}, vol.~12, no.~2, pp. 1224--1237, Mar. 2021.

\bibitem{4-1157}
M.~Nazemi, P.~Dehghanian, and M.~Lejeune, ``A mixed-integer distributionally robust chance-constrained model for optimal topology control in power grids with uncertain renewables,'' in \emph{2019 {IEEE} Milan {PowerTech}}.\hskip 1em plus 0.5em minus 0.4em\relax {IEEE}, Jun. 2019.

\bibitem{4-61}
B.~Kocuk, S.~S. Dey, and X.~A. Sun, ``New formulation and strong {MISOCP} relaxations for {AC} optimal transmission switching problem,'' \emph{IEEE Trans. Power Syst.}, vol.~32, no.~6, pp. 4161--4170, Nov. 2017.

\bibitem{4-1516}
T.~Lan, Z.~Zhou, W.~Wang, and G.~M. Huang, ``Stochastic optimization for {AC} optimal transmission switching with generalized benders decomposition,'' \emph{International Journal of Electrical Power {\&} Energy Systems}, vol. 133, p. 107140, dec 2021.

\bibitem{4-1517}
C.~Guo, H.~Nagarajan, and M.~Bodur, ``Tightening quadratic convex relaxations for the ac optimal transmission switching problem,'' Dec. 2022.

\bibitem{4-1518}
F.~Capitanescu, J.~M. Ramos, P.~Panciatici, D.~Kirschen, A.~M. Marcolini, L.~Platbrood, and L.~Wehenkel, ``State-of-the-art, challenges, and future trends in security constrained optimal power flow,'' \emph{Electric Power Systems Research}, vol.~81, no.~8, pp. 1731--1741, aug 2011.

\bibitem{4-1519}
M.~I. Alizadeh and F.~Capitanescu, ``A tractable linearization-based approximated solution methodology to stochastic multi-period {AC} security-constrained optimal power flow,'' \emph{{IEEE} Transactions on Power Systems}, pp. 1--13, 2022.

\bibitem{4-1520}
S.~Mhanna and P.~Mancarella, ``An exact sequential linear programming algorithm for the optimal power flow problem,'' \emph{{IEEE} Transactions on Power Systems}, vol.~37, no.~1, pp. 666--679, jan 2022.

\bibitem{4-490}
D.~K. Molzahn, I.~A. Hiskens \emph{et~al.}, ``A survey of relaxations and approximations of the power flow equations,'' \emph{Foundations and Trends{\textregistered} in Electric Energy Systems}, vol.~4, no. 1-2, pp. 1--221, 2019.

\bibitem{4-1213}
A.~Velloso, D.~Pozo, and A.~Street, ``Distributionally robust transmission expansion planning: A multi-scale uncertainty approach,'' \emph{IEEE Trans. Power Syst.}, vol.~35, no.~5, pp. 3353--3365, Sep. 2020.

\bibitem{4-1240}
C.~Coffrin, B.~Knueven, J.~Holzer, and M.~Vuffray, ``The impacts of convex piecewise linear cost formulations on {AC} optimal power flow,'' \emph{Electric Power Systems Research}, vol. 199, p. 107191, Oct. 2021.

\bibitem{4-1163}
C.~Coffrin and P.~V. Hentenryck, ``A linear-programming approximation of {AC} power flows,'' \emph{{INFORMS} Journal on Computing}, vol.~26, no.~4, pp. 718--734, Nov. 2014.

\bibitem{4-1129}
W.~E. Brown and E.~Moreno-Centeno, ``Transmission-line switching for load shed prevention via an accelerated linear programming approximation of {AC} power flows,'' \emph{IEEE Trans. Power Syst.}, vol.~35, no.~4, pp. 2575--2585, Jul. 2020.

\bibitem{4-1154}
J.~Ostrowski, J.~Wang, and C.~Liu, ``Exploiting symmetry in transmission lines for transmission switching,'' \emph{IEEE Trans. Power Syst.}, vol.~27, no.~3, pp. 1708--1709, Aug. 2012.

\bibitem{4-1239}
F.~Margot, ``Symmetry in integer linear programming,'' in \emph{50 Years of Integer Programming 1958-2008}.\hskip 1em plus 0.5em minus 0.4em\relax Springer Berlin Heidelberg, Nov. 2009, pp. 647--686.

\bibitem{4-1328}
R.~Chen, J.~Wang, A.~Botterud, and H.~Sun, ``Wind power providing flexible ramp product,'' \emph{IEEE Trans. Power Syst.}, vol.~32, no.~3, pp. 2049--2061, May 2017.

\bibitem{4-1330}
B.-I. Craciun, T.~Kerekes, D.~Sera, R.~Teodorescu, and U.~D. Annakkage, ``Power ramp limitation capabilities of large {PV} power plants with active power reserves,'' \emph{IEEE Trans. Sustain. Energy}, vol.~8, no.~2, pp. 573--581, Apr. 2017.

\bibitem{4-1355}
T.~Han, D.~J. Hill, and Y.~Song, ``Formulating connectedness in security-constrained optimal transmission switching problems,'' \emph{arXiv:2202.02805}, Feb. 2022.

\bibitem{4-1279}
M.~S. Bazaraa, H.~D. Sherali, and C.~M. Shetty, \emph{{Nonlinear Programming: Theory and Algorithms}}.\hskip 1em plus 0.5em minus 0.4em\relax John Wiley \& Sons, 2013.

\bibitem{4-1277}
G.~A. Hanasusanto, D.~Kuhn, and W.~Wiesemann, ``Computational complexity of stochastic programming problems,'' \emph{Mathematical Programming}, vol. 159, no. 1-2, pp. 557--569, Oct. 2015.

\bibitem{4-1322}
S.~Pineda and J.~M. Morales, ``Solving linear bilevel problems using {Big-Ms}: {Not} all that glitters is gold,'' \emph{IEEE Trans. Power Syst.}, vol.~34, no.~3, pp. 2469--2471, May 2019.

\bibitem{4-543}
S.~Fattahi, J.~Lavaei, and A.~Atamt{\"u}rk, ``A bound strengthening method for optimal transmission switching in power systems,'' \emph{IEEE Trans. Power Syst.}, vol.~34, no.~1, pp. 280--291, Jan. 2018.

\bibitem{4-1273}
\BIBentryALTinterwordspacing
T.~Kleinert and M.~Schmidt, ``Why there is no need to use a {Big-M} in linear bilevel optimization: {A} computational study of two ready-to-use approaches,'' 2020. [Online]. Available: \url{http://www.optimization-online.org/DB_FILE/2020/10/8065.pdf}
\BIBentrySTDinterwordspacing

\bibitem{4-971-1}
\BIBentryALTinterwordspacing
T.~Han, ``Structure-oriented optimization and control,'' 2021. [Online]. Available: \url{https://github.com/thanever/SOC/tree/master/Uncertainty/data}
\BIBentrySTDinterwordspacing

\bibitem{4-1356}
\BIBentryALTinterwordspacing
M.~e. Marko~Cepin, \emph{Assessment of Power System Reliability}.\hskip 1em plus 0.5em minus 0.4em\relax Springer London, Jul. 2011. [Online]. Available: \url{https://www.ebook.de/de/product/14756882/marko_cepin_marko_epin_assessment_of_power_system_reliability.html}
\BIBentrySTDinterwordspacing

\bibitem{4-1521}
\BIBentryALTinterwordspacing
GAMS, ``Gams/scenred2 documentation.'' [Online]. Available: \url{https://www.gams.com/42/docs/T_SCENRED2.html}
\BIBentrySTDinterwordspacing

\bibitem{4-1265}
\BIBentryALTinterwordspacing
L.~Zhao and B.~Zeng, ``An exact algorithm for two-stage robust optimization with mixed integer recourse problems,'' 2012. [Online]. Available: \url{https://optimization-online.org/2012/01/3310}
\BIBentrySTDinterwordspacing

\end{thebibliography}

\end{document}